\documentclass[14pt,a4paper,twoside]{amsart}
\usepackage{mathpazo}
\usepackage{amssymb}
\usepackage{mathtools}
\usepackage{parskip}
\usepackage{amsthm}
\usepackage{csquotes}
\usepackage{hyperref}
\usepackage{tikz}
\usepackage{tikz-cd}

\hypersetup{
colorlinks=true,
linkcolor=blue,
filecolor=blue,
urlcolor=blue,
citecolor=blue,
pdfpagemode=Fullscreen
}
\newtheorem{theorem}{Theorem}[section]
\newtheorem{lemma}[theorem]{Lemma}
\newtheorem*{Acknowledgement}{\textnormal{\textbf{Acknowledgement}}}
\theoremstyle{definition}
\newtheorem{definition}[theorem]{Definition}

\newtheorem{corollary}[theorem]{Corollary}
\newtheorem{proposition}[theorem]{Proposition}

\newtheorem{remark}[theorem]{Remark}

\numberwithin{equation}{section}

\newcommand{\beqa}{\begin{eqnarray*}}
\newcommand{\eeqa}{\end{eqnarray*}}
\newcommand{\beqn}{\begin{eqnarray}}
\newcommand{\eeqn}{\end{eqnarray}}

\newcounter{cnt1}
\newcounter{cnt2}
\newcounter{cnt3}
\newcommand{\blr}{\begin{list}{$($\roman{cnt1}$)$}
        {\usecounter{cnt1} \setlength{\topsep}{0pt}
                \setlength{\itemsep}{0pt}}}
\newcommand{\bla}{\begin{list}{$($\alph{cnt2}$)$}
        {\usecounter{cnt2} \setlength{\topsep}{0pt}
                \setlength{\itemsep}{0pt}}}
\newcommand{\bln}{\begin{list}{$($\arabic{cnt3}$)$}
        {\usecounter{cnt3} \setlength{\topsep}{0pt}
                \setlength{\itemsep}{0pt}}}
\newcommand{\el}{\end{list}}
\newtheorem{thm}{Theorem}

\newtheorem{Def}[thm]{Definition}

\newtheorem{rem}[thm]{Remark}
\newcommand{\Rem}{\begin{rem} \rm}
\newcommand{\bdfn}{\begin{Def} \rm}
\newcommand{\edfn}{\end{Def}}
			
\title{Monomial and binomial retracts of polynomial rings}

\author{Sagnik Chakraborty$^1*$}

\address{$^1$Department of Mathematics, 
	Ramakrishna Mission Vivekananda Educational and Research Institute, 
	Belur Math,  Howrah 711202,
	West Bengal, India}		 

\email{ jusagnik28@gmail.com}

\author{Madhuparna Pal$^2$}
\address {$^2$Department of Mathematics,
	Ramakrishna Mission Vivekananda Educational and Research Institute, 
	Belur Math,  Howrah 711202,
	West Bengal, India}

\email{honey.madhu777p@gmail.com}

\begin{document}
	\maketitle
	
	\iffalse
	\author{Sagnik Chakraborty$*$}\\
	
	\address{Department of Mathematics, 
		Ramakrishna Mission Vivekananda Educational and Research Institute, 
		Belur Math,  Howrah 711202,
		West Bengal, India}		 
	
	\email{Email: jusagnik28@gmail.com}
	
	\author{Madhuparna Pal}\\
	\address {Department of Mathematics,
		Ramakrishna Mission Vivekananda Educational and Research Institute, 
		Belur Math,  Howrah 711202,
		West Bengal, India}
	
	\email{Email: honey.madhu777p@gmail.com}	
\fi	
	
\footnote{$*$ corresponding author}

\begin{abstract}
Let $R$ be a ring and $B \coloneqq R[X_1, \dots, X_n]$ the polynomial ring in $n$ variables over $R$. In this article, we consider retractions $\varphi : B \longrightarrow B$ such that $\varphi(X_i)$ is $0$, or a monomial or the sum of two monomials. We prove that, under certain conditions on the base ring $R$, the resulting retracts are also polynomial rings over $R$. We characterize different monomial retractions of $B$, i.e., where $\varphi(X_i)$ is $0$ or a monomial for all $i$, giving the same retract. We also generalize the structure of monomial retracts by allowing a few variables to map to arbitrary polynomials, but only under some additional restrictions on the base ring $R$.
\end{abstract}
 \smallskip
 
 \noindent
 {\small {{\bf Keywords}. Polynomial ring, Monomial, Binomial, Retraction, Symmetric algebra, $\mathbb A^n$-fibration.}
 	
 }

\section{Introduction}
Throughout this article, all rings are assumed to be commutative rings with unity. If $R$ is a ring, we use $R^{[n]}$ to denote the polynomial ring in $n$ variables over $R$.  
If $A \subseteq B$ are $R$-algebras, then $A$ is said to be an \textit{$R$-algebra retract} (or an {\it $R$-retract}) of $B$ if one (and hence all) of the following equivalent conditions is satisfied.
\begin{enumerate}
\item There exists an $R$-algebra homomorphism $ \varphi : B \longrightarrow B$ such that $\varphi^2 = \varphi$ and $\varphi(B)=A$; such a map $\varphi$ is called an \textit{$R$-algebra retraction}, or an {\it $R$-retraction}, for short.
	\item There exists an $R$-algebra homomorphism $ \varphi : B \longrightarrow A$ such that $\varphi|_A = \text{id}_A$.
	\item There exists an ideal $I$ of $B$ such that $B= A \oplus I$ as $A$-modules.
\end{enumerate}
Recall that a {\it monomial} of $B=R^{[n]}:=R[X_1,\dots,X_n]$ is a polynomial of the form $\lambda X_1^{a_1}\dots X_n^{a_n}$, where $\lambda\in R$ is a nonzero element and $a_1,\dots,a_n$ are non-negative integers; and such a monomial is called {\it monic} if $\lambda=1$. Note that, according to this definition, a nonzero constant polynomial is a monomial. An $R$-retraction $\varphi:B\longrightarrow B$ is said to be a ({\it monic}) \textit{monomial $R$-retraction} if $\varphi(X_i)$ is $0$ or a (monic) monomial for all $i$, and the image of a (monic) monomial $R$-retraction is called a ({\it monic}) {\it monomial $R$-retract}.\\
Similarly, an $R$-retraction $\varphi:B\longrightarrow B$ is called a {\it binomial $R$-retraction} if $\varphi(X_i)$ is a sum of at most two monomials for all $i$; and the resulting $R$-retract is called a {\it binomial $R$-retract} (see \ref{binomialretract} for details).
%Note that monomial $R$-retractions are special cases of binomial $R$-retractions.\\

\par In \cite{costa1977retracts} D. L. Costa raised the following question. 
\begin{center}
	{\it Is every $k$-retract of $k^{[n]}$ a polynomial ring over $k$ where $k$ is a field?}		
\end{center}
Costa himself proved in \cite[Theorem 3.4, 3.5]{costa1977retracts} that if $k$ is a field and $n\le 2$, then all $k$-retracts of $k^{[n]}$ are polynomial rings over $k$. The answer to the above question is negative in positive characteristic for all $n\ge 4$, as shown by Gupta in \cite{gupta2014cancellation}, \cite{gupta2014zariski}. In characteristic zero, it was independently proved by Nagamine (\cite[Theorem 2.5]{nagamine2019note}), and by the first author, Dasgupta, Dutta and Gupta (\cite[Theorem 5.8]{chakraborty2021some}) that all $k$-retracts of $k^{[3]}$ are polynomial rings over $k$. However, the question remains open for $n=3$ in positive characteristic, and for all $n\ge 4$ in characteristic zero.

\par In \cite{kojima2024remarks} Kojima, Nagamine and Sasagawa gave some sufficient conditions for a $k$-retract of $k^{[3]}$ to be a polynomial ring over $k$ in positive characteristic (Lemma 4.3, Proposition 4.4, 4.6), and described all monic monomial $k$-retractions of $k^{[n]}$ for $n=2,3$, which turn out to be polynomial rings over $k$ (Example 4.2). Gupta and Nagamine proved in \cite[Theorem 1.2]{gupta2023retracts} that if $R$ is an integral domain, then every $R$-retract of the Laurent polynomial ring $R^{[\pm n]}:=R[X_1^{\pm 1},\dots,X_n^{\pm 1}]$ is also a Laurent polynomial ring over $R$ for all $n\ge 0$. Their proof relies on the fact that every subgroup of a finitely generated free abelian group $G$ is also a finitely generated free abelian group, which is a direct summand of $G$. However, in the case of polynomial rings, the exponents of the monomials form a finitely generated free monoid, and we know that a submonoid of a finitely generated free monoid need not be free. Therefore, the same idea does not work for polynomial rings.

\par The above works suggest that although Costa's question about $k$-retracts of $k^{[n]}$ being polynomial rings or not remains open, in general, in characteristic zero, and has a negative answer, in general, in positive characteristic, it is worth exploring the nature of such retracts for some special classes of retractions, like monomial or binomial retractions, of polynomial rings over $k$, where $k$ is a field or, more generally, even an arbitrary ring.

\par In this paper, we investigate the structure of monomial and binomial $R$-retracts of polynomial rings over $R$. We prove that 
\begin{enumerate}
    \item {\it If $R$ does not contain any nontrivial idempotent elements, then a monomial $R$-retract of $R^{[n]}$ is a polynomial ring over $R$} (Theorem \ref{monthm1}).
    \item {\it If $R$ is a UFD, then a binomial $R$-retract of $R^{[n]}$ is a polynomial ring over $R$} (Corollary \ref{binomcor}).
\end{enumerate}

In fact, our results describe the structure of monomial and binomial $R$-retracts of $R^{[n]}$ in terms of symmetric algebras of projective $R$-modules over more general rings $R$ than the ones mentioned above. We also characterize different monomial $R$-retractions that result in the same retract (Theorem \ref{diffthm}). One can slightly relax the condition for monomial $R$-retractions by allowing a few variables to map to arbitrary polynomials. In this paper, we have studied $R$-retractions of $R^{[n]}$ where all but at most three variables map to zero or monomials. We give a complete description of such retracts, although under some additional hypothesis on the base ring $R$ (Proposition \ref{genprop1}, and Corollaries \ref{diffcor} and \ref{gencor1}). 

We will now discuss the layout of the paper. In Section \ref{sec:prelim}, we set up the notation and recall some definitions and known results. In Section \ref{sec:mon}, we study the structure of monomial $R$-retractions of $R^{[n]}$. In Section \ref{sec:diff}, we classify different monomial $R$-retractions of $R^{[n]}$ that result in the same retract. In Section \ref{sec:gen}, we generalize our results for the structure of monomial $R$-retractions by allowing a few variables to map to arbitrary polynomials, albeit at the cost of certain restrictions on the base ring $R$. In Section \ref{sec:binom}, we describe the structure of binomial $R$-retracts of $R^{[n]}$ for an integral domain $R$, and in Section \ref{sec:open}, we discuss some open problems.

\section{Preliminaries}
\label{sec:prelim}

{\bf Definitions and Notation:}

If $A\subseteq B$ are rings, then $b_1,\ldots,b_n\in B$ are said to be {\it algebraically independent} (over $A$) if the $A$-algebra homomorphism $A[X_1,\ldots,X_n]\longrightarrow B$, which sends $X_i$ to $b_i$ for all $i$, is injective.\\
Recall that if $f\in R[X_1,\ldots,X_n]$ is a nonzero polynomial, then the {\it degree} (respectively, {\it order}) of $f$, denoted by $\text{ deg }f$ (respectively, $\text{ord }f$), is defined as the maximum (respectively, the minimum) of the total degrees of monomials present in $f$. Note that if $f$ is a monomial, then $\text{deg }f=\text{ord }f$.\\% the total degree 
Let $A\subseteq B$ be integral domains. Then the {\it transcendence degree} of the field of fractions of $B$ over that of $A$ is denoted by $\text{tr.deg}_A B$.\\
Let $R$ be a ring. A finitely generated flat $R$-algebra $A$ is said to be an {\it $ \mathbb A^n$-fibration} over $ R $ if, for each prime ideal $\mathfrak p $ of $ R $, $ A \otimes_R \kappa(\mathfrak p) = \kappa(\mathfrak p)^{[n]} $, where $\kappa(\mathfrak p)$ is the residue field of $R$ at the prime ideal $\mathfrak p$. An $\mathbb A^n$-fibration $A$ over $R$ is called a {\it trivial} $\mathbb A^n$-fibration if $A\cong R^{[n]}$ as $R$-algebras.\\
We use $\mathcal M_n(R)$ to denote the ring of $n\times n$ matrices over $R$; and $E\in \mathcal M_n(R)$ is called an {\it idempotent} matrix if $E^2=E$. Recall that, in general, if $R$ is a ring then $e\in R$ is called an {\it idempotent} element if $e^2=e$; and an idempotent element $e$ is called a {\it non-trivial idempotent element} if $e\neq 0,1$.\\
If $R$ is a ring and $M$ an $R$-module, then $\text{Sym}_R(M)$ denotes the symmetric algebra of $M$ over $R$.\\
A {\it (positively) graded ring} $R$ is a ring $R$ together with a direct sum decomposition of $R$ as an additive group $R=\bigoplus_{i\in \mathbb N}R_i$ such that $R_iR_j\subseteq R_{i+j}$ for all non-negative integers $i,j$. A non-zero element $r\in R$ is said to be {\it homogeneous} if $r\in R_i$ for some $i\in \mathbb N$, and in this case, the non-negative integer $i$ is called the {\it degree} of $r$. The ideal of $R$ generated by the homogeneous elements of positive degree is called
the {\it irrelevant ideal} of $R$ and is denoted by $R_+$, so that $R=R_0\oplus R_+$.\\
A ring $R$ is said to be a {\it semi-normal} ring if, for all $b,c\in R$ that satisfy $b^2=c^3$, there exists an element $a\in R$ such that $b=a^3$ and $c=a^2$. An integral domain that is also semi-normal is called a {\it semi-normal domain}.\\
If $R$ is a ring, then the {\it Picard group} of $R$, denoted by $\mathrm{Pic }\ R$, is defined as the (abelian) group of isomorphism classes of projective $R$-modules of rank one, where the binary operation is the tensor product of $R$-modules. Note that every projective module of rank one is finitely generated.

{\bf Preliminary results:}

The following two results were proved by Costa in \cite{costa1977retracts}. The first result characterizes retracts of polynomial rings over a UFD in some special cases (\cite[Theorem 3.5 and the subsequent remark]{costa1977retracts}), and the second describes the structure of $R$-algebra retracts of the polynomial ring $R[X]$ (\cite[Theorem 3.4]{costa1977retracts}).

\begin{theorem}
	\label{co1}
	Let $R$ be a UFD and let $A$ be a retract of $ B \coloneqq R[x_1, \dots, x_n] (=R^{[n]})$.  
    \begin{enumerate}
        \item If $ \text{tr.deg}_R A =0$, then $A=R$.
        \item If $\text{tr.deg}_RA=1$, then $A=R^{[1]}$.
        \item If $\text{tr.deg}_RA=n$, then $A=B$.
    \end{enumerate}

    \end{theorem}

    \begin{theorem}
        \label{co2}
        Let $R$ be a ring and $R[X]$ the polynomial ring in one variable over $R$. If $\varphi:R[X]\longrightarrow R[X]$ is an $R$-algebra retract, then $\varphi(X)=a+eX$ for some $a,e\in R$ such that $e$ is an idempotent element, and $ae=0$. In particular, every $R$-algebra retract of $R[X]$ is of the form $R[eX]$, where $e$ is an idempotent element of $R$.
    \end{theorem}

    Next, we state a well-known result about idempotent elements of polynomial rings. We also include a proof here for the lack of ready reference. However, it should be noted that similar and more general results were proved by Gilmer in \cite[Section 10]{gilmer}.

    \begin{proposition}
        \label{gilmer}
        {\it If $A:=A_0\oplus A_1\oplus\ldots$ is a positively graded ring, then every idempotent element of $A$ is contained in $A_0$. In particular, if $A$ is a polynomial ring over $R$, then all idempotent elements of $A$ are contained in $R$.}
    \end{proposition}

    \begin{proof}
        If possible, let $a\in A\setminus A_0$ be an idempotent element with homogeneous decomposition $a=a_0+a_1+\ldots+a_n$, where $a_n\neq 0$. Let $r$ be the smallest positive integer such that $a_r\neq 0$. From $a^2=a$, we get $a_0^2=a_0$; and comparing the $r$-th homogeneous component of the two sides, we see that $a_r=2a_0a_r$. Since $2a_0a_r=2a_0^2a_r$, it follows that $a_r=a_0a_r=2a_0a_r$. Therefore, $a_r=a_0a_r=0$, which is a contradiction.
    \end{proof}

The following standard properties of symmetric algebras can be found in \cite[Propositions 7 and 9 in section 6 of Chapter III]{baki}.

\begin{proposition}
	\label{baki}
    {\it Let $R$ be a ring.
    \begin{enumerate}
        \item If $M$ is an $R$-module and $A$ is an $R$-algebra, then $A\otimes_R\text{Sym}_R(M)\cong \text{Sym}_A(A\otimes_RM)$ as $A$-algebras. 
        \item If $M_1,\ldots, M_n$ are $R$-modules, then \[\text{Sym}_R(M_1\oplus\cdots\oplus M_n)\cong \text{Sym}_R(M_1)\otimes_R\cdots\otimes_R\text{Sym}_R(M_n)\] as $R$-algebras.
    \end{enumerate}}

 Swan proved the following result for Picard groups of polynomial rings in \cite[Theorem 1]{swan}.

    \begin{theorem}
        \label{swan}
        If $R$ is a semi-normal ring, then the natural map $\mathrm{Pic }\ R\longrightarrow \mathrm{Pic }\ R[X_1,\ldots,X_n]$, obtained by extension of scalars, is an isomorphism of groups for all $n\ge 1$. 
    \end{theorem}

\end{proposition}

The following five consecutive results are taken from \cite{chakraborty2021some}.  

\begin{theorem}
	\label{5.4}
    \begin{enumerate}
        \item Let $R\subseteq A$ be integral domains, with $R$ being Noetherian. If $ \text{tr.deg}_R A =1$, then $A$ is an $ \mathbb A^1 $-fibration over $ R $ if and only if it is an $ R$-algebra retract of $ R^{[m]} $ for some positive integer $m$ (\cite[Theorem 5.4]{chakraborty2021some}).
        \item Let $R $ be a Noetherian  semi-normal domain, $B=R^{[n]}$ and $A$ a retract of $B$ such that $\text{tr.deg}_R A=1$. Then $A \cong \text{Sym}_R(J)$ for some invertible ideal $J$ of $R$ (\cite[Corollary 5.6]{chakraborty2021some}). This was previously proved by C. Greither (\cite[Theorem 2.3]{greither}).
        \item Let $ R $ be a Noetherian domain containing $ \mathbb Q $, and $ A $ an integral domain containing $ R $ such that $ \text{tr.deg}_R A =2 $. Then $ A $ is an $ \mathbb A^2 $-fibration over $ R $ if and only if it is an $ R $-algebra retract of $ R^{[m]} $ for some positive integer $ m $ (\cite[Theorem 5.9]{chakraborty2021some}).
        \item Let $ R $ be an integral domain and let $ A $ be a graded $ R $-subalgebra of $ B\coloneqq R[X_1, \dots, X_n] $ with standard grading. Suppose that there exists a retraction $ \varphi : B \longrightarrow B $ such that $ \pi(B_+) \subseteq B_+ $. If $A$ is the image of $\varphi$, then $ A = \text{Sym}_R(M) $ for some finitely generated projective $ R $-submodule $ M$ of the free $R$-module $ RX_1 \oplus \dots \oplus RX_n $ (\cite[Theorem 5.13]{chakraborty2021some}).
        \item Let $ A $ be a subring of $ B $ and $ \varphi : B \longrightarrow B$ a retraction whose image is $A$. Then for any ideal $J \subseteq \text{ker } \varphi $ of $ B $, $ A $ is the image of the induced retraction $\bar\varphi: B/J \longrightarrow B/J$. In particular, $A$ is a retract of $B/J$ (\cite[Lemma 6.1(ii)]{chakraborty2021some}).
        
    \end{enumerate}

\end{theorem}

\section{Monomial Retractions of Polynomial Rings}
\label{sec:mon}

Let $n$ be a positive integer, $R$ a ring and $B \coloneqq R[X_{1}, \dots , X_n]$ the polynomial ring in $n$ variables over $R$. In general, a monomial $R$-retract of $B$ may not be a polynomial ring. In \cite[Theorem 3.4]{costa1977retracts} it was proved that {\it an $R$-algebra retract of $R[X]$ is of the form $R[eX]$, where $e$ is an idempotent element of $R$}. Clearly, if $e$ is a non-trivial idempotent element of $R$, then $R[eX]$ is not isomorphic to a polynomial ring over $R$. In general, if $e\in R$ is a non-trivial idempotent element, then it is easy to construct monomial $R$-retracts of $B$ that are not polynomial rings over $R$. We will show that if $R$ does not contain any non-trivial idempotent elements, then every monomial $R$-retract of $B$ is a polynomial ring over $R$.\\

We first make the following simple observation about monomial retractions.

\begin{lemma}
    \label{monlem0}
    Let $R$ be a ring and let $B:=R[X_1,\ldots,X_n]$ be the polynomial ring in $n$ variables over $R$. Let $\varphi :B\longrightarrow B$ be a monomial $R$-retraction. Suppose that $\varphi(X_j)=af(\neq 0)$ for some $j$, where $f$ is a monic monomial that is divisible by $X_j$. Then the following assertions hold.
    \begin{enumerate}
        \item The monomial $f$ is not divisible by $X_j^2$, and every variable other than $X_j$ that appears in $f$ is mapped to a nonzero element of $R$ under $\varphi$.
        \item The ideal $(a)$ is an idempotent ideal; in other words, there exists a nonzero idempotent element $e\in R$ such that $(a)=(e)$.
    \end{enumerate}
     
\end{lemma}

\begin{proof}
    Let $f=\prod_{i=1}^nX_i^{e_i}$, where each $e_i$ is a non-negative integer and $e_j\ge1$. Since $\varphi(X_i)$ is a monomial for all $i$, and $af=\varphi(X_j)=\varphi^2(X_j)=a\big(\prod_{i=1}^n\varphi(X_i)^{e_i}\big)=a^{e_j+1}f^{e_j}\big(\prod_{i\neq j}\varphi(X_i)^{e_i}\big)$, comparing the degrees of the two sides, we see that $e_j=1$, and if $e_i\ge 1$ for some $i\neq j$, then $\varphi(X_i)$ is a nonzero element of $R$. Finally, comparing the coefficients of the two sides, we see that $a=a^2r$, where $r:=\prod_{i\neq j, e_i\ge 1}\varphi(X_i)^{e_i}$. Clearly, $e:=ar$ is an idempotent element and $(a)=(e)$.
\end{proof}

If $\varphi$ is a monomial $R$-retraction of the polynomial ring $R[X_1,\ldots,X_n]$, then the variables $X_i$ for which $\varphi(X_i)$ is a nonzero multiple of $X_i$ play an important role in the description of the associated retract. So, we make the following definition.

\begin{definition}
Let $R$ be a ring and $B:=R[X_1,\ldots,X_n]$ be the polynomial ring in $n$ variables over $R$. Let $\varphi:B\longrightarrow B$ be a monomial $R$-retraction. We define the {\it support of $\varphi$}, denoted by $\mathcal S_\varphi$, as
\[\mathcal S_\varphi:=\{i\in \{1,\ldots,n\}\ |\ \varphi(X_i) \text{ is a nonzero monomial that is divisible by }X_i\}.\]
Note that if $i\in \{1,\ldots,n\}\setminus\mathcal S_\varphi$, then either $X_i\in \text{ker }\varphi$, or $\varphi(X_i)$ is a nonzero element of $R$ or $\varphi(X_i)$ is a non-constant monomial in which $X_i$ does not appear. In view of {\it lemma \ref{monlem0}}, if $i\in\mathcal S_\varphi$ then $\varphi(X_i)$ is of the form $X_i\cdot M_i$, where $M_i$ is a nonzero monomial, and if $X_j$ appears in $M_i$, for some $j\neq i$, then $\varphi(X_j)$ is a nonzero element of $R$.

\end{definition}

Next, we prove the following special case of our result for monomial retractions, which has a similar flavor to Theorem \ref{5.4}\textcolor{blue}{(iv)}.

\begin{proposition}    
   \label{monprop1}
    {\it Let $R$ be a ring and $B:=R[X_1,\ldots,X_n]$ the polynomial ring in $n$ variables over $R$. Let $\varphi :B\longrightarrow B$ be a monomial $R$-retraction with image $A$. If $\varphi(X_i)\in B_+$, i.e., if $\varphi(X_i)=0$ or $\varphi(X_i)$ is a non-constant monomial, for all $i$, then $\mathcal S_\varphi=\{i_1,\ldots,i_m\}\subseteq \{1,\ldots,n\}$, for some $0\le m\le n$, such that $\varphi(X_{i_j})=\lambda_{i_j}X_{i_j}$ for some nonzero idempotent element $\lambda_{i_j}\in R$ for all $j=1,\ldots,m$, and $A=R[\lambda_{i_1}X_{i_1},\ldots,\lambda_{i_m}X_{i_m}]$.}
\end{proposition}

\begin{proof}
Let $\varphi_i:=\varphi(X_i)$ for all $i=1,\ldots, n$. If $\varphi_i\neq0$ for some $i$, then we can write $\varphi_i$ as $\phi_i=\lambda_if_i$, where $\lambda_i$ is a nonzero element of $R$ and $f_i$ is a monic monomial. Clearly, each variable that appears in $f_i$ must map to a linear monomial; for otherwise, the degree of $\varphi(f_i)$ would be greater than the degree of $f_i$, leading to a contradiction. Let $\Omega_\varphi:=\{i\in \{1,\ldots,n\}\ | \ \varphi_i\text{ is a linear monomial}\}$. If $i\in \Omega_\varphi$ and $\varphi_i=\lambda_iX_{i'}$ for some $i'\in \{1,\ldots,n\}$, then $\lambda_iX_{i'}=\varphi(\lambda_iX_{i'})=(\lambda_i\cdot\lambda_{i'})X_{i'}$, where $\varphi_{i'}=\lambda_{i'}X_{i'}$, implying that $i'\in \Omega_\varphi$. Therefore, we obtain a set-theoretic function $\phi:\Omega_\varphi\to\Omega_\varphi$ such that $\varphi_i$ is a nonzero scalar multiple of $X_{\phi(i)}$ for all $i\in \Omega_\varphi$. As $\varphi^2=\varphi$, it follows that $\phi^2=\phi$. Then $\mathcal S_\varphi=\text{im }\phi$ is the set of fixed points of $\phi$, i.e., $\mathcal S_\varphi=\{i\in \Omega_\varphi\ |\ \varphi_i\text{ is a nonzero scalar multiple of }X_i\}$. For each $i\in \mathcal S_\varphi$, $\varphi_i=\lambda_iX_i=\varphi(\lambda_iX_i)=\lambda_i^2X_i$, which implies that $\lambda_i$ is a nonzero idempotent element of $R$. If $\mathcal S_\varphi=\{i_1,\ldots,i_m\}$, then $0\le m\le n$, and it is clear from the discussion that $A$, as an $R$-algebra, is generated by $\{\varphi_{i_j}\}_{j=1}^m$, i.e., $A=R[\lambda_{i_1}X_{i_1},\ldots,\lambda_{i_m}X_{i_m}]$.

\end{proof}

\begin{corollary}
    \label{moncor1}
    Let $R$ be a ring and $B:=R[X_1,\ldots,X_n]$ the polynomial ring in $n$ variables over $R$. Let $\varphi:B\longrightarrow B$ be a monomial $R$-retraction with image $A$ such that $\varphi(X_i)$ is a nonzero monomial that is divisible by $X_i$ for all $i=1,\ldots,n$. Then there exist nonzero idempotent elements $\lambda_1,\ldots,\lambda_n\in R$ such that $\varphi(X_i)=\lambda_iX_i$ for all $i=1,\ldots,n$. In particular, $A=R[\lambda_1X_1,\ldots,\lambda_nX_n]$.  
\end{corollary}

\begin{proof}
    The assertion follows immediately from {\it lemma \ref{monlem0}} and {\it proposition \ref{monprop1}}.
\end{proof}

\begin{remark}
    \label{monrem0}
    With the set-up and notation as in {\it proposition \ref{monprop1}}, note that the retract $A$ is a monomial $R$-subalgebra of $B$. It follows from the proof of {\it proposition \ref{monprop1}} that for each $i\in \{1,\ldots,n\}\setminus\mathcal S_\varphi$, if $\varphi_i\neq 0$, then $\varphi_i$ is a monomial in $\lambda_{i_1}X_{i_1},\ldots,\lambda_{i_m}X_{i_m}$, which is also a consequence of the standard properties of monomial algebras.
\end{remark}

\begin{theorem}
    \label{monthm1}
    Let $R$ be a ring and $B:=R[X_1,\ldots,X_n]$ the polynomial ring in $n$ variables over $R$. Let $\varphi :B\longrightarrow B$ be a monomial $R$-retraction with image $A$. Let $\mathcal S_\varphi=\{i_1,\ldots,i_m\}\subseteq \{1,\ldots,n\}$, for some $0\le m\le n$. Then $A=R[e_{1}f_{1},\ldots,e_{m}f_{m}]$ for some nonzero idempotent elements $e_{1},\ldots,e_{m}\in R$ such that the following assertions hold.    
    
    \begin{enumerate}
    \item The polynomials $f_{1},\ldots,f_{m}$ are monic monomials that are algebraically independent over $R$.
        
        \item For each $j=1,\ldots, m$, $\varphi(X_{i_j})=\lambda_jf_j$ such that $(\lambda_j)=(e_j)$ is a nonzero idempotent ideal, so that $A=R[\varphi(X_{i_1}),\ldots,\varphi(X_{i_m})]$; and each $f_j$ is of the form $X_{i_j}\cdot M_j$, where $M_j$ is a monomial and every variable that appears in $M_j$ is mapped to a nonzero element of $R$ under $\varphi$.
        \item If $i\in \{1,\ldots, n\}\setminus\mathcal S_\varphi$, then either $\varphi(X_i)\in R$ or it is a nonzero monomial in $e_1f_1,\ldots,e_mf_m$.       
        
    \end{enumerate}
    
    In particular, if $R$ does not contain any non-trivial idempotent elements, or if $\varphi$ is a monic monomial $R$-retraction, then $A=R^{[m]}$, for some $0\le m\le n$.
    
\end{theorem}

\begin{proof}
   We apply induction on $n$, the number of variables. The case $n=1$ is a special case of {\it theorem \ref{co2}}. Suppose that the statement is true for some positive integer $r\ge 1$, and let $n=r+1$. If $\varphi(X_i)\in B_+$ for all $i$, then the assertion follows from {\it proposition \ref{monprop1}}. In fact, in this case, the monomials $M_j$'s can be chosen from among the variables $X_i$'s and the coefficients $\lambda_j$'s are nonzero idempotent elements. So, after re-indexing the variables, if necessary, we may assume that $c:=\varphi(X_{r+1})$ is a nonzero element of $R$. Let $J\subseteq B$ be the ideal generated by $X_{r+1}-c$. Since $J\subseteq \text{ker }\varphi$, by {\it theorem \ref{5.4}\textcolor{blue}{(v)}}, we get an induced $R$-retraction $\bar\varphi:B/J\longrightarrow B/J$ whose image, namely $\bar A:=\bar\varphi(B/J)$, is isomorphic to $A$ as an $R$-algebra; in fact, the natural projection $\pi:B\longrightarrow B/J$ induces an $R$-algebra isomorphism between $A$ and $\bar A$. Note that $B/J\cong R[X_1,\ldots,X_r]=R^{[r]}$, and $\bar\varphi$ is a monomial $R$-retraction.\\ 
   Let $\varphi_i:=\varphi(X_i)$ and $\bar\varphi_i:=\bar\varphi(X_i)$ for all $i=1,\ldots,r$. Note that $\bar\varphi_i$ is obtained by plugging $c$ in place of $X_{r+1}$ in $\varphi_i$ for all $i=1,\ldots,r$. By the induction hypothesis, there exists a non-negative integer $m$ such that $\mathcal S_{\bar\varphi}=\{i_1,\ldots,i_m\}\subseteq \{1,\ldots,r\}$ and  $\bar A=R[\bar\varphi_{i_1},\ldots,\bar\varphi_{i_m}]$, where each $\bar\varphi_{i_j}$ is a nonzero monomial whose coefficient generates a nonzero idempotent ideal of $R$. Let $\varphi_{i_j}=\lambda_jf_j$, where $\lambda_j\in R$ and $f_j$ is a monic monomial. Then $\bar\varphi_{i_j}=(\lambda_jc^{\epsilon_j})g_j$, where $\epsilon_j$ is the largest non-negative integer such that $X_{r+1}^{\epsilon_j}$ divides $f_j$, and $g_j$ is obtained by plugging $1$ in place of $X_{r+1}$ in $f_j$. By the induction hypothesis, $g_1,\ldots,g_m$ are algebraically independent over $R$. Therefore, $f_1,\ldots,f_m$ are also algebraically independent over $R$. Since $X_{i_j}$ divides $g_j$ for all $j$, it follows that $X_{i_j}$ also divides $f_j$ for all $j$, implying that $\mathcal S_{\bar\varphi}\subseteq \mathcal S_\varphi$. Conversely, suppose that $i\in \mathcal S_\varphi$. Then $X_i$ divides $\varphi(X_i)$, and, therefore, also divides $\bar \varphi(X_i)$. However, in order to conclude that $i\in \mathcal S_{\bar\varphi}$, we must show that $\bar\varphi_i\neq 0$. Let $\varphi_i:=\lambda\prod_{l=1}^{r+1}X_l^{e_{i,l}}$ for some nonzero $\lambda\in R$, where $e_{i,l}$ are non-negative integers such that $e_{i,i}\ge 1$. If $\bar\varphi_i=\lambda c^{i,r+1}\prod_{l=1}^rX_l^{e_{i,l}}=0$, then $\lambda c^{e_{i,r+1}}=0$. But then $\varphi_i=\varphi(\varphi_i)=\lambda c^{e_{i,r+1}}\prod_{l=1}^r\varphi(X_l)^{e_{i,l}}=0$, which is a contradiction. Therefore, we conclude that $\mathcal S_\varphi=\mathcal S_{\bar\varphi}$. Since $X_{i_j}$ divides $\varphi_{i_j}$ for all $j=1,\ldots ,m$, it follows from {\it lemma \ref{monlem0}} that each $\lambda_j$ generates a nonzero idempotent ideal of $R$, and each $f_j$ is of the form $X_{i_j}\cdot M_j$, where $M_j$ is a monomial such that every variable appearing in $M_j$ is mapped to a nonzero element of $R$ under $\varphi$. We claim that $A=R[\varphi_{i_1},\ldots,\varphi_{i_m}]$. Since $\bar A$ is generated by the $\bar\varphi_{i_j}$'s, if $f\in A$ then we can find a polynomial $f_0$ is $m$ variables such that $\pi(f)=f_0\big(\bar\varphi_{i_1},\ldots,\bar\varphi_{i_m}\big)$, implying that $f-f_0\big(\varphi_{i_1}\ldots,\varphi_{i_m}\big)\in J$. Since $A\cap J=(0)$, $f=f_0(\varphi_{i_1},\ldots,\varphi_{i_m})\in R[\varphi_{i_1},\ldots,\varphi_{i_m}]$, and consequently, $A=R[\varphi_{i_1},\ldots,\varphi_{i_m}]=R[\lambda_1f_1,\ldots,\lambda_mf_m]$. By {\it lemma \ref{monlem0}}, for each $j=1,\ldots,m$, there exists a nonzero idempotent element $e_j\in R$ such that $(\lambda_j)=(e_j)$, which implies that $A=R[e_1f_1,\ldots,e_mf_m]$. Since $A$ is a monomial $R$-subalgebra of $B$, the third assertion follows from the standard properties of monomial algebras. \\   
   The last assertion is obvious because if $R$ does not contain any non-trivial idempotent elements, then each $\lambda_j$ must be a unit, and if $\varphi$ is a monic monomial $R$-retraction, then $\lambda_j=1$ for all $j=1,\ldots,m$. In these two cases, $A=R[\varphi(X_{i_1}),\ldots,\varphi(X_{i_m})]\cong R^{[m]}$, where $\mathcal S_\varphi=\{i_1,\ldots,i_m\}$.

   \end{proof}

   \begin{remark}
       \label{monrem1}
       Let $R$ be a ring, $B:=R^{[n]}$ and $\varphi:B\longrightarrow B$ a monomial $R$-retraction. 
       \begin{enumerate}
           \item If the irrelevant ideal of $B$ is not invariant under $\varphi$, then the retract $A:=\varphi(B)$ may not be of the form described in {\it proposition \ref{monprop1}}. For example, if $n=2$, $\varphi(X_1):=X_1X_2$ and $\varphi(X_2):=1$, then $A=R[X_1X_2]$.
           \item In {\it theorem \ref{monthm1}}, since $f_1,\ldots,f_m$ are algebraically independent over $R$, it follows that the retract $A=R[e_1f_1,\ldots,e_mf_m]$ is isomorphic to the symmetric $R$-algebra $\text{Sym}_R(Re_1\oplus\cdots\oplus Re_m)$.
           \item The constants $\lambda_j$'s, as described in {\it theorem \ref{monthm1}}, may not themselves be idempotent elements. For example, let $R$ be a ring that contains a non-trivial idempotent element $e$ and a unit $u$ such that $e(u-1)\neq 0$. Let $\varphi:R[X_1,X_2]\longrightarrow R[X_1,X_2]$ be the monomial $R$-retraction, defined as $\varphi(X_1):=euX_1X_2$ and $\varphi(X_2):=u^{-1}$. Then $eu$ itself is not an idempotent element, but $(eu)=(e)$ is an idempotent ideal.
           \item In {\it theorem \ref{monthm1}}, if $R$ is an integral domain, then $A=R[\varphi(X_i)]_{i\in \mathcal S_\varphi}\cong R^{[m]}$, where $m=|\mathcal S_\varphi|$. Let $J\subseteq B$ be the ideal generated by the set $\{X_i-\varphi(X_i)\}_{i\in \{1,\ldots,n\}\setminus\mathcal S_\varphi}$. Since $J\subseteq\text{ker }\varphi$ and $B/J$ is naturally isomorphic to $R[X_i]_{i\in \mathcal S_\varphi}\cong R^{[m]}$, by {\it theorem \ref{5.4}}\textcolor{blue}{(v)}, we obtain an induced monomial $R$-retraction $\bar\varphi$ of $R^{[m]}$, such that $A$, as an $R$-algebra, is isomorphic to $\bar A:=\bar\varphi\big(R^{[m]}\big)$. In particular, $\text{tr.deg}_R\bar A=m$. It is well-known that $\bar A$ is algebraically closed in $R^{[m]}$ (for  example, see \cite[Lemma 1.3]{costa1977retracts}). Therefore, $\bar A=R^{[m]}$ and $\bar\varphi$ is the identity map of $R^{[m]}$; and consequently, $\text{ker }\varphi=J$.
       \end{enumerate}

\end{remark}       

\section{Different Monomial Retractions with the Same Retract}
\label{sec:diff}

In this section, we discuss how different monomial $R$-retractions of a polynomial ring over $R$ result in the same retract. We maintain the same set-up and notation as in {\it section \ref{sec:mon}}, i.e., $R$ is a ring, $B:=R[X_1,\ldots,X_n]$ is a polynomial ring in $n$ variables over $R$ and $\varphi:B\longrightarrow B$ is a monomial $R$-retraction with image $A$. For any such $\varphi$, we define $\mathcal S_\varphi$ to be the set of positive integers $i\in \{1,\ldots,n\}$ such that $\varphi(X_i)$ is a nonzero multiple of $X_i$.\\ 

Since we are now mainly interested in the retracts, there is no harm in getting rid of the `extra' variables. Note that if a variable, say $X_n$, does not appear in $A$, then $B_n:=R[X_1,\ldots,X_{n-1}]\subseteq B$ is $\varphi$-invariant, and it follows from {\it theorem \ref{monthm1}} that the induced monomial $R$-retraction $\varphi|_{B_n}:B_n\longrightarrow B_n$ has the same image as $\varphi$. On the other hand, we can extend $\varphi$ to a monomial $R$-retraction $\tilde \varphi :B[Y_1,\ldots,Y_m]\to B[Y_1,\ldots,Y_m]$ by sending each $Y_i$ to zero or a monomial of $B$, which is contained in $A$. Then $\varphi$ and $\tilde\varphi$ have the same image. However, no $Y_i$ appears in the image of $\tilde\varphi$. We can thus continue adding more and more variables to $B$ without changing the retract. Therefore, to study different monomial retractions with the same retract, we may ignore all variables that do not appear in $A$.%do this by artificially adding as many variables as we wish. On the other hand, 

\begin{definition}
    Let $R$ be a ring and $B:=R[X_1,\ldots,X_n]$ the polynomial ring in $n$ variables over $R$. Let $\varphi:B\longrightarrow B$ be a monomial $R$-retraction. We call $\varphi$ {\it non-degenerate} if each $X_i$ appears in the image of $\varphi$, or equivalently, if each $X_i$ appears in $\varphi(X_j)$ for some $j$. A monomial $R$-retraction that is not non-degenerate, is called {\it degenerate}.
\end{definition}

In the rest of this section, we will only consider monomial $R$-retractions that are non-degenerate. Let $\varphi:B\longrightarrow B$ be a non-degenerate monomial $R$-retraction, where $B:=R[X_1,\ldots,X_n]$. If $\varphi(X_i)=0$ for some $i$, then $X_i$ cannot appear in $\varphi(X_j)$ for all $j$. Since $\varphi$ is non-degenerate, we can, therefore, assume that $\varphi(X_i)\neq 0$ for all $i$. By {\it theorem \ref{monthm1}}, $A=R[\lambda_iX_iM_i]_{i\in \mathcal S_\varphi}$, where $\varphi(X_i):=\lambda_iX_iM_i$ for all $i\in \mathcal S_\varphi$, each $\lambda_i$ is a nonzero element of $R$ and each $M_i$ is a monic monomial in variables that are not contained in $\mathcal S_\varphi$. Note that if $c_i:=\varphi(M_i)$ for some $i\in \mathcal S_\varphi$, then $\lambda_i=\lambda_i^2c_i$.\\ 
If each $\lambda_i$ is a unit, then $A=R[X_iM_i]\cong R^{[m]}$, where $m:=|\mathcal S_\varphi|$. In this case, for any family of invertible elements $\tilde \mu :=\{\mu_i\}_{i\notin\mathcal S_\varphi}$ of $R$, we can define a monomial $R$-retraction $\varphi_{\tilde\mu}:B\longrightarrow B$ as
\begin{equation*}
\varphi_{\tilde \mu}(X_i)=
    \begin{cases}
        \mu_i  & \text{if } i\notin \mathcal S_\varphi\\
        \frac{X_iM_i}{\varphi_{\tilde\mu}(M_i)} & \text{ otherwise. } %x \in \mathbb{R}\setminus\mathbb{Q}
    \end{cases}
\end{equation*}
Then it is easy to check that $\varphi_{\tilde \mu}(B)=A$. Therefore, if $\mathcal S_{\varphi}$ is a proper subset of $\{1,\ldots,n\}$ and $R$ contains infinitely many units, then for any family of non-constant monic monomials $\{M_i\}_{i\in \mathcal S_\varphi}$, there exist infinitely many monomial $R$-retractions of $B$ with the image $A=R[X_iM_i]_{i\in \mathcal S_\varphi}$.\\

However, the situation improves drastically if we restrict ourselves to monic monomial $R$-retractions. So, let $\varphi:B\longrightarrow B$ be a monic monomial $R$-retraction. We want to find the set of monic monomial $R$-retractions of $B$ that have the same image as $\varphi$, say $A$. Let $\psi:B\longrightarrow B$ be another monic monomial $R$-retraction with the image $A$. Since $\varphi,\psi$ are monic monomial $R$-retractions, by {\it theorem \ref{monthm1}}, $A\cong R^{[m]}$ as $R$-algebras, where $m:=|\mathcal S_\varphi|=|\mathcal S_\psi|$. We know that $A=R[\varphi(X_i)]_{i\in \mathcal S_\varphi}=R[\psi(X_j)]_{j\in \mathcal S_\psi}$. Let $\mathbf {M}_\varphi:= \{\varphi(X_i)\}_{i\in \mathbf S_\varphi}$ and $\mathbf M_\psi:=\{\psi(X_j)\}_{j\in \mathcal S_\psi}$. We claim that $\mathbf M_\varphi=\mathbf M_\psi$.\\ 
Let $\varphi(X_i):=X_iM_i$ and $\psi(X_j):=X_jM'_j$ for all $i\in \mathcal S_\varphi$ and $j\in \mathcal S_\psi$. Then each $M_i$ (respectively, $M'_j$) is a monic monomial in variables that are not contained in $\mathcal S_\varphi$ (respectively, not contained in $\mathcal S_\psi$). Since $R[\varphi(X_i)]_{i\in \mathcal S_\varphi}=R[\psi(X_j)]_{j\in \mathcal S_\psi}$, every monomial in $\mathbf M_\varphi$ can be written as a monic monomial in the elements of $\mathbf M_\psi$, and vice versa. Let $\mathbf{M}_\varphi\cap\mathbf M_\psi$ contain $m_0$ elements, with $0\le m_0\le m$. If $m_0=m$, then we are done. Otherwise, let $\varphi(X_i)$ and $\psi(Xj)$ be two monomials of the least degree among the monomials of $\mathbf M_\varphi\setminus(\mathbf M_\varphi\cap\mathbf M_\psi)$ and $\mathbf M_\psi\setminus(\mathbf M_\varphi\cap\mathbf M_\psi)$, respectively. Note that $\varphi(X_i)$ is not divisible by $\varphi(X_{i'})$ for all $i'\in \mathcal S_\varphi\setminus\{i\}$. Similarly, $\psi(X_j)$ is not divisible by $\psi(X_{j'})$ for all $j'\in \mathcal S_\psi\setminus\{j\}$. Therefore, we can write $\varphi(X_i)$ and $\psi(X_j)$ as a monomial in the elements of $\mathbf M_\psi\setminus(\mathbf M_\varphi\cap\mathbf{M}_\psi)$ and $\mathbf M_\varphi\setminus(\mathbf M_\varphi\cap\mathbf M_\psi)$, respectively. Consequently, $\text{deg }\varphi(X_i)=\text{deg }\psi(X_j)$, and it follows that $\varphi(X_i)\in\mathbf M_\psi$ (and $\psi(X_j)\in \mathbf M_\varphi$), which is a contradiction. Therefore, $\mathbf M_\varphi=\mathbf M_\psi$, which results in a set-theoretic bijection $\sigma:\mathbf M_\varphi\to \mathbf M_\psi$ such that $\varphi(X_i)=\psi\big(X_{\sigma(i)}\big)$ for all $i\in \mathbf M_\varphi$. Since $\varphi(X_i)=X_iM_i$, either $\sigma(i)=i$ or $X_{\sigma(i)}$ appears in $M_i$. Since $X_{\sigma(i)}$ cannot appear in $\psi(X_j)$ for all $j\neq \sigma(i)$, we see that, in the latter case, $X_{\sigma(i)}$ is a factor of $M_i$ and $X_{\sigma(i)}$ does not appear in $M_{i'}$ for all $i'\neq i$. Therefore, for each $i\in \mathcal S_\varphi$, if we define $\Gamma_i$ as the set of all positive integers $s\in \{1,\ldots,n\}$ such that $X_s$ appears in $\varphi(X_i)$, but does not appear in $\varphi(X_{i'})$ for all $i'\neq i$, then $\sigma(i)\in \Gamma_i$ for all $i\in \mathcal S_\varphi$. Note that $i\in \Gamma_i$ for all $i\in \mathcal S_\varphi$.\\
Conversely, it is easy to see that if $\sigma:\mathcal S_\varphi\to \{1,\ldots,n\}$ is an injective map such that $\sigma(i)\in \Gamma_i$ for all $i\in \mathcal S_\varphi$, then the monic monomial $R$-retraction $\varphi_\sigma:B\longrightarrow B$, defined as

\begin{equation*}
\varphi_{\sigma}(X_j)=
    \begin{cases}
        \varphi(X_i)  & \text{if } j=\sigma(i)\text{ for some }i\in \mathcal S_\varphi\\
        1 & \text{ otherwise, } %x \in \mathbb{R}\setminus\mathbb{Q}
    \end{cases}
\end{equation*}
has the same image as $\varphi$.\\

We record our observations in the form of the following theorem.

\begin{theorem}
    \label{diffthm}
    Let $R$ be a ring and $\varphi,\psi$ non-degenerate monic monomial $R$-retractions of the polynomial ring $B:=R[X_1,\ldots,X_n]$. Then $\varphi$ and $\psi$ have the same image iff there exists a set-theoretic bijection $\sigma:\mathcal S_\varphi\to \mathcal S_\psi$ such that $\varphi(X_i)=\psi\big(X_{\sigma(i)}\big)$ for all $i\in \mathcal S_\varphi$. In this case, $\sigma(i)\in \Gamma_i$ for all $i\in \mathcal S_\varphi$, where $\Gamma_i$ is the set of all $s\in \{1,\ldots,n\}$ such that $X_s$ appears in $\varphi(X_i)$, but $X_s$ does not appear in $\varphi(X_{i'})$ for all $i'\neq i$.\\
    In particular, the number of distinct monic monomial $R$-retractions of $B$ that have the same image as $\varphi$ is given by $\prod_{i\in \mathcal S_\varphi}|\Gamma_i|$.
\end{theorem}

The above result can be used to count the number of distinct monic monomial retractions with the same image, as well as to construct monic monomial retractions with the same retract. For example, if $m,n$ are positive integers, then the image of the monic monomial $R$-retraction $\varphi:R[X_1,X_2,X_3]\longrightarrow R[X_1,X_2,X_3]$, given by $\varphi(X_1):=X_1X_3^m$, $\varphi(X_2):=X_2X_3^n$ and $\varphi(X_3):=1$, cannot be obtained as the image of a different monic monomial $R$-retraction.\\
On the other hand, the two monic monomial $R$-retractions $\varphi,\psi:R[X_1,X_2,X_3]\longrightarrow R[X_1,X_2,X_3]$, defined as $\varphi(X_1):=X_1X_3$, $\varphi(X_2):=X_2$ and $\varphi(X_3):=1$, and $\psi(X_1):=1$, $\psi(X_2):=X_2$ and $\psi(X_3):=X_1X_3$ have the same image, namely $R[X_1X_3,X_2]$.

\section{Miscellaneous Generalizations}
\label{sec:gen}

In this section, we generalize some of the results that were proved for monomial retractions by allowing a few variables to map to arbitrary nonzero polynomials, not necessarily monomials. However, generalizations often come at a cost, as we will eventually place certain restrictions on the base ring $R$.\\

As usual, let $R$ be a ring and $B:=R[X_1,\ldots,X_n]$ the polynomial ring in $n$ variables over $R$. Let $\varphi :B\longrightarrow B$ be an $R$-retraction with image $A$. Our first result is a partial generalization of {\it theorem \ref{monthm1}} where we show that if $R$ does not contain any non-trivial idempotent elements, then $A$ is isomorphic a polynomial ring over $R$ even if one of the variables is mapped to an arbitrary polynomial, which is not zero or a monomial. Then we will further generalize our result by proving that for an arbitrary ring $R$, such a retract $A$ is isomorphic to the symmetric algebra of an $R$-module $M$, which is a finite direct sum of cyclic projective $R$-modules (Corollary \ref{diffcor}).

\begin{proposition}
\label{genprop1}
	{\it Let $R$ be a ring that does not contain any non-trivial idempotent elements. Let $B:=R[X_1,\ldots,X_n]$ be the polynomial in $n$ variables over $R$, and $\varphi: B\longrightarrow B$ an $R$-retraction with image $A$. If $\varphi(X_i)$ is either zero or a monomial for all but at most one $i$, then $A\cong R^{[m]}$, as $R$-algebras, for some $0\le m\le n$.}
    
\end{proposition}

\begin{proof}
	We apply induction on $n$, the number of variables. If $n=1$ and $\varphi(X)\in R$, then $A=R$; otherwise, by {\it theorem \ref{co2}}, $A$ is of the form $R[eX]$ for a nonzero idempotent element $e\in R$. Since $R$ does not contain any non-trivial idempotent elements, it follows that $e=1$ and $A=R[X]$. Next, we assume that the assertion is true for some positive integer $r\ge 1$, and let $n=r+1$. If $\varphi(X_i)$ is zero or a monomial for all $i=1,\ldots,r+1$, then there is nothing to prove, as the assertion follows from {\it theorem \ref{monthm1}}. So, without loss of generality, we may assume that $\varphi(X_i)$ is zero or a monomial for all $i=1,\ldots,r$, and $\varphi(X_{r+1})$ is a nonzero polynomial which is not a monomial.\\ 
    If $X_i$ does not appear in $\varphi(X_i)$ for some $i\le r$, then let $J_i\subseteq B$ be the principal ideal generated by $X_i-\varphi(X_i)$. As $J_i\subseteq\text{ker }\varphi$, using {\it theorem \ref{5.4}\textcolor{blue}{(v)}}, we can consider the induced retraction $\bar\varphi:B/J_i\longrightarrow B/J_i$. The natural projection $B\longrightarrow B/J_i$ induces an $R$-algebra isomorphism between $A$ and $\bar A$, the image of $\bar\varphi$. Note that $B/J_i$ is naturally isomorphic to the polynomial ring $R[X_1,\ldots,X_{i-1},X_{i+1},\ldots,X_{r+1}]\cong R^{[r]}$. Since $\bar\varphi:R^{[r]}\longrightarrow R^{[r]}$ maps every variable other than possibly $X_{r+1}$ to zero or a monomial, by the induction hypothesis, $A\cong \bar A\cong R^{[m]}$, for some $0\le m\le r$; therefore, the claim is true.\\
    So, we assume that $\varphi(X_i)$ is a nonzero monomial that is divisible by $X_i$ for all $i=1,\ldots,r$. Then each $\varphi(X_i)$ can be written as $\varphi(X_i)=(\lambda_iX_i)\cdot M_i$, where $\lambda_i$ is a nonzero element of $R$ and $M_i$ is a monic monomial. Since $\varphi(X_i)=\varphi^2(X_i)$, $(\lambda_iX_i)M_i=(\lambda_i^2X_iM_i)\cdot(\varphi(M_i))$; and with $X_iM_i$ being a non-zero divisor, it follows that $\lambda_i=\lambda_i^2\varphi(M_i)$. If $X_j$ appears in $M_i$ for some $j\le r$, then the order of $\varphi(M_i)$ would be at least one, which is a contradiction. Therefore, $M_i=X_{r+1}^{\epsilon_i}$ for some non-negative integer $\epsilon_i$. 
    Since $\varphi$ is a retraction and $\varphi(X_i)=\lambda_iX_iX_{r+1}^{\epsilon_i}$, we get $\lambda_iX_iX_{r+1}^{\epsilon_i}=\varphi(\lambda_iX_iX_{r+1}^{\epsilon_i})=\lambda_i^2X_iX_{r+1}^{\epsilon_i}\varphi(X_{r+1})^{\epsilon_i}$. With $X_iX_{r+1}^{\epsilon_i}$ being a non-zero divisor, it follows that $\lambda_i=\lambda_i^2\varphi(X_{r+1})^{\epsilon_i}$. Since $R$ does not contain any non-trivial idempotent elements and  $\lambda_i\varphi(X_{r+1})^{\epsilon_i}$ is a nonzero idempotent element of $R$, we conclude that $\lambda_i\varphi(X_{r+1})^{\epsilon_i}=1$ for all $i=1,\ldots,r$. In particular, if $\epsilon_i=0$ for some $i$, then $\lambda_i=1$.\\ 
    If $\epsilon_i=0$ for all $i=1,\ldots,r$, then $\varphi$ fixes each $X_i$. As a result, $\varphi$ can be realized as an $R[X_1,\ldots,X_r]$-retraction of $B$. Therefore, by {\it theorem \ref{co2}}, $\varphi(Y)=a+eY$ for some $a,e\in R[X_1,\ldots,X_r]$, such that $e$ is an idempotent element, and $ae=0$. Since $R$
    does not contain any non-trivial idempotent elements, by {\it proposition \ref{gilmer}}, the same is true for $R[X_1,\ldots,X_r]$, which implies $e=1$ or $e=0$. In the first case, $\varphi(Y)=Y$, implying that $\varphi=\text{id}_B$; and in the second case, $A=R'=R[X_1,\ldots,X_r]=R^{[r]}$.\\
    Next, we assume that $\epsilon_j>0$ for some $j$. Let $\varphi(X_{r+1})=c+g$, where $c\in R$ and $g\in B_+$. From $\varphi(X_j)=\varphi^2(X_j)$, we get 
    \begin{align}\lambda_jX_jX_{r+1}^{\epsilon_j}=\lambda_j^2X_jX_{r+1}^{\epsilon_j}(c+g)^{\epsilon_j}\ . \end{align}
    Since $X_jX_{r+1}^{\epsilon_j}$ is a non-zero divisor, $\lambda_j=\lambda_j^2(c+g)^{\epsilon_j}$, which implies that $\lambda_j(c+g)^{\epsilon_j}$ is a nonzero idempotent element of $B$. Since $R$ does not contain any non-trivial idempotent elements, by {\it proposition \ref{gilmer}}, the same is true for $B$. Therefore, $\lambda_j(c+g)^{\epsilon_j}=1$. It follows that $\lambda_jc^{\epsilon_j}=1$, and that each coefficient of $g$ is a nilpotent element of $R$. Let $R':=R[X_1X_{r+1}^{\epsilon_1},\ldots,X_rX_{r+1}^{\epsilon_r}]\subseteq A$. Let $R[T_1,\ldots,T_r]$ be a polynomial ring in $r$ variables over $R$, and let $\psi:R[T_1,\ldots,T_r]\longrightarrow B$ be the $R$-algebra homomorphism that sends $T_i$ to $X_iX_{r+1}^{\epsilon_i}$ for all $1\le i\le r$. If $\pi:B\longrightarrow R[X_1,\ldots,X_r]$ is the $R[X_1,\ldots,X_r]$-algebra homomorphism, which sends $X_{r+1}$ to $1$, then $\pi\circ \psi$ is an $R$-algebra isomorphism, which implies that $\psi:R[T_1,\ldots,T_r]\longrightarrow R'$ is also an isomorphism.\\
    We claim that $\varphi(X_{r+1})=c+g\in R'$, or equivalently, $A=R'\cong R^{[r]}$. Let $\mathcal N_g\subseteq R$ be the ideal that is generated by the coefficients of $g$. Then $\mathcal N_g$ is a nilpotent ideal, i.e., $(\mathcal N_g)^s=0$ for all $s>>0$. Let $M:=aX_1^{e_1}\cdot\ldots\cdot X_{r}^{e_{r}}X_{r+1}^{e_{r+1}}$ be a nonzero monomial of $B$. Then 
    \[\varphi(M)=a(\lambda_1X_1X_{r+1}^{\epsilon_1})^{e_1}\cdot\ldots\cdot (\lambda_rX_rX_{r+1}^{\epsilon_r})^{e_r}(c+g)^{e_{r+1}}=\big(ac^{e_{r+1}}\prod_{i=1}^{r}\lambda_i^{e_i}\big)X_1^{e_1}\cdot\ldots\cdot X_r^{e_r}X_{r+1}^{\sum_{i=1}^r\epsilon_ie_i}+P_M,\]
    where $P_M$ is a polynomial whose all coefficients are contained in $a\mathcal N_g$. If $M_0:=\varphi(M)-P_M$,
    %$M_0:=\big(ac^{e_{r+1}}\prod_{i=1}^{r}\lambda_i^{e_i}\big)X_1^{e_1}\cdot\ldots\cdot X_r^{e_r}X_{r+1}^{\sum_{i=1}^rd_ie_i}$, 
    then it is clear that $M_0\in R'$. Since every polynomial can be (uniquely) written as a sum of monomials, if $h\in B$, then $\varphi(h)$ can be written as $\varphi(h)=h_0+P_h$ such that $h_0\in R'$ and $P_h$ is a polynomial whose every coefficient is contained in $\mathcal J_h\mathcal N_g$, where $\mathcal J_h\subseteq R$ is the ideal generated by the coefficients of $h$. Note that $B=R'\oplus R''$ as $R$-modules, where $R''$ is the free $R$-submodule of $B$ generated by the set of monomials $\{X_1^{l_1}\cdot\ldots\cdot X_r^{l_r}X_{r+1}^{l_{r+1}}\ |\ l_{r+1}\neq\sum_{i=1}^r\epsilon_il_i\}$. As $R'\subseteq A$, $A=R'\oplus(A\cap R'')$ as $R$-modules. Therefore, $A=R'$ iff $A\cap R''=0$. If possible, let $A\cap R''\neq 0$. Then we can choose a nonzero polynomial $f\in A\cap R''$. We can write $f=\varphi(f)$ uniquely as the sum of two polynomials $f_0+P_f$, where $f_0\in R'$ and $P_f\in R''$. However, since $f\in R''$, $f_0=0$, or equivalently, $f=P_f$. If $\mathcal J_f\subseteq B$ is the ideal generated by the coefficients of $f$, then it follows from the above discussion that every coefficient of $P_f$ is contained in $\mathcal J_f\mathcal N_g$. Since $f=P_f$, it follows that $\mathcal J_f=\mathcal J_f\mathcal N_g$. With $\mathcal N_g$ being a nilpotent ideal, it implies that $\mathcal J_f=0$, or equivalently, $f=0$. Therefore, $f\in R'$, and we conclude that $A=R'\cong R^{[r]}$.
\end{proof}

Next, we generalize the above result to an arbitrary ring $R$ that may contain nontrivial idempotent elements. For that, we need the following elementary result about symmetric algebras. 

\begin{lemma}
    \label{genelem}
    Let $R_1,\ldots,R_n$ be rings, and let $R:=R_1\times\cdots\times R_n$. Let $M_i$ be an $R_i$-module for all $i=1,\ldots,n$. We can view each $M_i$ as an $R$-module by restriction of scalars. Then
    \[\text{Sym}_R(M_1\oplus\cdots\oplus M_n)\cong\prod_{i=1}^n\text{Sym}_{R_i}(M_i)\]
    as $R$-algebras.
\end{lemma}

\begin{proof}
    Let $e_i\in R$ be the idempotent element that has $1$ in the $i$-th co-ordinate and zero everywhere else. For each $i$, let $\mathcal J_i\subseteq R$ be the ideal generated by the set $\{e_j\}_{j\neq i}$. Then $R_i=Re_i=R/\mathcal J_i$ for all $i=1,\ldots,n$. Since $\mathcal J_i+\mathcal J_j=R$ for all $i\neq j$, it follows that if $N$ is an $R_i$-module for some $i$, then $R_j\otimes_RN=0$ for all $j\neq i$. Now, by {\it proposition \ref{baki}},
    \[\text{Sym}_R(M_1\oplus\cdots\oplus M_n)\cong R\otimes_R\text{Sym}_R(M_1\oplus\cdots\oplus M_n)\cong (R_1\times\cdots\times R_n)\otimes_R\text{Sym}_R(M_1\oplus\cdots\oplus M_n)\]
    \[\cong\prod_{i=1}^nR_i\otimes_R\big(\text{Sym}_R(M_1\oplus\cdots\oplus M_n)\big)\cong \prod_{i=1}^n\text{Sym}_{R_i}\big(R_i\otimes_R(M_1\oplus\cdots\oplus M_n)\big)\cong\prod_{i=1}^n\text{Sym}_{R_i}(M_i), \]
    as $R_i\otimes_RM_j=0$ for all $j\neq i$, and $R_i\otimes_RM_i\cong M_i$.
\end{proof}

\begin{corollary}
    \label{diffcor}
    Let $R$ be a ring and $B:=R[X_1,\ldots,X_n]$ the polynomial ring in $n$ variables over $R$. Let $\varphi:B\longrightarrow B$ be an $R$-retraction with image $A$. If $\varphi(X_i)$ is zero or a monomial for all but at most $i$, then there exists an $R$-module $M$, which is the direct sum of $m$ cyclic projective $R$-modules, for some non-negative integer $m\le n$, such that $A\cong \text{Sym}_R(M)$ as $R$-algebras.
\end{corollary}

\begin{proof}
    Let $R'$ be a Noetherian subring of $R$ that contains the coefficients of $\varphi(X_i)$ for all $i=1,\ldots,n$, and let $B':=R'[X_1,\ldots,X_n]$. Then $\varphi$ induces an $R'$-retraction $\varphi':B'\longrightarrow B'$ such that $\varphi'(X_i)=\varphi(X_i)$ is zero or a monomial for all but at most one $i$. Let $A'$ be the image of $\varphi'$. As $B=R\otimes_{R'}B'$, $B=A\oplus \text{ker }\varphi$ and $B'=A'\oplus\text{ker }\varphi'$, it follows that $\text{ker }\varphi=R\otimes_{R'}\text{ker }\varphi'$ and $A=R\otimes_{R'}A'$.\\    
    Since $R'$ is a Noetherian ring, it contains finitely many minimal prime ideals, whose intersection is $\text{Nil }R'$, the nil radical of $R'$. It follows that $R'/\text{Nil }R'$ can be embedded in a finite product of integral domains, which implies that $R'/\text{Nil }R'$ contains only finitely many idempotent elements. Now, it is well-known that the natural projection $R'\to R'/\text{Nil }R'$ induces a bijection between the set of idempotent elements of $R'$ and the set of idempotent elements of $R'/\text{Nil }R'$ (for example, see \cite[Proposition 7.14]{Jacobson}). Consequently, the number of idempotent elements of $R'$ is also finite. So, we can choose nonzero idempotent elements $e_1,\ldots,e_s\in R'$, such that
    \begin{enumerate}%[(a)]
        \item $e_1+\cdots+e_s=1$,
        \item $e_ie_j=0$ for all $i\neq j$, and
        \item no $e_i$ can be written as the sum of two nonzero idempotent elements of $R'$.
    \end{enumerate}
    For each $i$, let $\mathcal J'_i\subseteq R'$ be the ideal generated by the set $\{e_j\}_{j\neq i}$, and let $\mathcal J_i:=R\mathcal J'_i$. Let $R'_i:=R'e_i\cong R'/\mathcal J'_i$ and $R_i:=Re_i\cong R/\mathcal J_i$ for all $i=1,\ldots,s$. Then $R'=R'_1\times \cdots\times R'_s$ and $R=R_1\times \cdots\times R_s$. For each $i$, $\varphi$ (respectively, $\varphi'$) induces an $R_i$-retraction $\varphi_i:B_i\longrightarrow B_i$ (respectively, an $R'_i$-retraction $\varphi'_i:B_i'\longrightarrow B'_i$), where $B_i:=R_i[X_1,\ldots,X_n]$ (respectively, $B'_i:=R'_i[X_1,\ldots,X_n]$). Let $A_i$ (respectively, $A'_i$) be the image of $\varphi_i$ (respectively, $\varphi'_i$) for all $i=1,\ldots,s$. Then $A=A_1\times\cdots\times A_s$ and $A'=A'_1\times\cdots\times A'_s$. It is also clear that $A_i=R_i\otimes_{R'_i}A'_i$ for all $i=1,\ldots,s$. Note that each $\varphi_i'$ has the same property as $\varphi'$, namely, $\varphi'_i(X_j)$ is zero or a monomial for all but at most one $j$. Since $R'_i$ does not contain any non-trivial idempotent elements, by {\it proposition \ref{genprop1}}, $A'_i$ is isomorphic to a polynomial ring over $R'_i$ in $d_i$ variables, for some non-negative integer $d_i\le n$. Since $A_i\cong R_i\otimes_{R'_i}A_i$, it follows that $A_i\cong R_i^{[d_i]}=\text{Sym}_{R_i}(R_i^{d_i})$. We re-index the idempotent elements, if necessary, to assume that $d_1\le d_2\le\ldots\le d_s$. Then it follows from {\it lemma \ref{genelem}} that
    \[A=\prod_{i=1}^sA_i\cong \prod_{i=1}^s\text{Sym}_{R_i}(R_i^{d_i})\cong\text{Sym}_R(R_1^{d_1}\oplus\cdots\oplus R_s^{d_s}).\]
    For each $i=1\ldots,s$, let $\tilde e_i:=\sum_{j=i}^se_i$. Then every $\tilde e_i$ is an idempotent element of $R$, and $R\tilde e_i=\prod_{j=i}^sR_i$ is a cyclic projective $R$-module. If $P_i:=R\tilde e_i$ for all $i=1,\ldots,s$, then %. Then each $P_i$ is a cyclic projective $R$-module, and
    \[A\cong \text{Sym}_R\big(R_1^{d_1}\oplus\cdots\oplus R_s^{d_s}\big)=\text{Sym}_R\big(P_1^{d_1}\oplus P_2^{d_2-d_1}\oplus P_3^{d_3-d_2}\oplus\cdots\oplus P_s^{d_s-d_{s-1}}\big).\]
    Let $d_0:=0$ and $m:=d_s$. Then it is clear that $0\le m\le n$ and $M:=\bigoplus_{i=1}^s P_i^{d_i-d_{i-1}}$ is the direct sum of $m$ cyclic projective $R$-modules, and $A\cong \text{Sym}_R(M)$.   
    
\end{proof}

Now we explore the structure of $R$-retracts of polynomial rings where all but at most two or three variables map to zero or monomials. For this purpose, in the remainder of the section, we assume that $R$ is an integral domain.

\begin{lemma}
    \label{genlem1}
    Let $R$ be an integral domain and $B:=R[X_1,\ldots,X_n,Y_1,\ldots,Y_m]$ the polynomial ring in $m+n$ variables over $R$. Let $\varphi: B\longrightarrow B$ be an $R$-retraction such that $\varphi(X_i)$ is zero or a monomial for all $i=1,\ldots,n$, and $\varphi(Y_j)$ is neither zero nor a monomial for all $j=1,\ldots,m$. If $B_0:=R[X_1,\ldots,X_n]$, then $B_0$ is $\varphi$-invariant, i.e., $\varphi(B_0)\subseteq B_0$. 
\end{lemma}

\begin{proof}
It is sufficient to show that $Y_j$ does not appear in $\varphi(X_i)$ for all $i=1,\ldots,n$ and $j=1,\ldots,m$. Let $i\in \{1.\ldots,n\}$. Since $R$ is an integral domain and $\varphi(X_i)$ is a monomial, every factor of $\varphi(X_i)$ is also a monomial. If $Y_j$ appears in $\varphi(X_i)$ for some $j$, then $\varphi(Y_j)$ divides $\varphi(X_i)$. Therefore, since $\varphi(Y_j)$ is not a monomial, $Y_j$ cannot appear in $\varphi(X_i)$ for all $j=1,\ldots,m$. Consequently, $\varphi(X_i)\in B_0$ for all $i=1,\ldots,n$, or equivalently, $\varphi(B_0)\subseteq B_0$.
\end{proof}

\begin{remark}
    With the same set-up and notation as in {\it lemma \ref{genlem1}}, with the exception of the assumption that $R$ is an integral domain, we give a few examples to show that one cannot, in general, eliminate the condition that $R$ is an integral domain.
    \begin{enumerate}
        \item Let $R$ be a ring with a non-trivial idempotent element $e\in R$. Let $\varphi:R[X,Y]\longrightarrow R[X,Y]$ be the $R$-retraction given by $\varphi(X):=eXY$ and $\varphi(Y):=e+(1-e)Y$. Then $B_0:=R[X]$ is not $\varphi$-invariant. Here, the retract $A=R[eXY,(1-e)Y]$ is not isomorphic to a polynomial ring over $R$; it is isomorphic to the symmetric $R$-algebra $\text{Sym}_R\big(Re\oplus R(1-e)\big)$.
        \item Let $R$ be a ring and $a\in R$ a nonzero element such that $a^2=0$. Let \[\varphi:R[X_1,X_2,Y_1,Y_2]\to R[X_1,X_2,Y_1,Y_2]\] be the $R$-retraction, given by $\varphi(X_1):=X_1$, $\varphi(X_2):=X_2Y_1Y_2$, $\varphi(Y_1):=1+a X_1$ and $\varphi(Y_2):=1-aX_1$. Then $B_0:=R[X_1,X_2]$ is not $\varphi$-invariant. Here, the retract $A=R[X_1,X_2Y_1Y_2]$ is isomorphic to $R^{[2]}$.
        \item Let $R$ be a ring with characteristic $p>0$, where $p$ is a prime number, and $a\in R$ a nonzero element that satisfies $a^p=0$. Let $\varphi:R[X_1,X_2,Y]\longrightarrow R[X_1,X_2,Y]$ be the $R$-retraction, defined as $\varphi(X_1):=X_1$, $\varphi(X_2):=X_2Y^p$ and $\varphi(Y):=1+aX_1$. Then $B_0:=R[X_1,X_2]$ is not $\varphi$-invariant. Here, the retract $A=R[X_1,X_2Y^p]$ is isomorphic to $R^{[2]}$.
        \item Let $R$ be a ring and $a,b\in R$ nonzero elements such that $ab=0$. Let $\varphi:R[X_1,X_2,Y]\longrightarrow R[X_1,X_2,Y]$ be the $R$-retraction, defined as $\varphi(X_1):=0$, $\varphi(X_2):=aY$ and $\varphi(Y):=Y+bX_1$. Then $B_0:=R[X_1,X_2]$ is not $\varphi$-invariant. Here, the retract $A=R[Y+bX_1]$ is isomorphic to $R^{[1]}$.
    \end{enumerate}
\end{remark}

The basic idea that we use to generalize the structure of monomial retracts is to `suppress' the variables that map to zero or a monomial by replacing the base ring $R$ with a polynomial ring over it. For that, we introduce the following definition.

\begin{definition}
    Let $R$ be an integral domain and $B:=R[X_1,\ldots,X_n]$ the polynomial ring in $n$ variables over $R$. Let  $\varphi: B\longrightarrow B$ be an $R$-retraction with image $A$. We define \[\Lambda_\varphi:=\{i\in\{1,\ldots,n\}\ |\ \varphi(X_i)\text{ is zero or a monomial}\}.\]
\end{definition}

\begin{proposition}
\label{genprop2}
{\it Let $R$ be an integral domain and $B:=R[X_1,\ldots,X_n]$ the polynomial ring in $n$ variables over $R$. Let $\varphi: B\longrightarrow B$ be an $R$-retraction with image $A$. If $\Lambda_\varphi$ contains $d$ elements, then there exist non-negative integers $m_0,n_0$, with $m_0\le n-d$ and $n_0\le d$, such that $A$, as an $R$-algebra, is isomorphic to a $B_0$-retract of the polynomial ring $B_0^{[m_0]}$, where $B_0:=R^{[n_0]}$.}
   
\end{proposition}

\begin{proof}
    We apply induction on $n$, the number of variables. Note that if $\Lambda_\varphi=\emptyset$, or equivalently, if $d=0$, then there is nothing to prove, as we can take $m_0:=n$ and $n_0:=0$. So, we assume, at the outset, that $d\ge 1$. If $n=d=1$, then there are two possibilities, either $\varphi(X_1)\in R$ or $\varphi(X_1)$ is a non-constant monomial. In the first case, $A=R$, so we can take $m_0=n_0:=0$ and $B_0:=R$; and if $\varphi(X_1)$ is a non-constant monomial, then it follows from {\it theorem \ref{co2}} that $\varphi$ is the identity map on $B$, so we can take $m_0:=0$ and $n_0:=1$, which implies $B_0:=B$. Now suppose that the claim is true for some positive integer $r\ge 1$, and let $n=r+1$. If $\varphi(X_i)\in R$ for some $i\in \Lambda_\varphi$, then let $J_i\subseteq B$ be the ideal generated by $X_i-\varphi(X_i)$. As $J_i\subseteq\text{ker }\varphi$, by {\it theorem \ref{5.4}\textcolor{blue}{(v)}}, we get an induced $R$-retraction $\bar\varphi:B/J_i\longrightarrow B/J_i$, and the natural projection $B\longrightarrow B/J_i$ induces an $R$-algebra isomorphism between $A$ and $\bar A:=\bar\varphi(B/J_i)$. Now $B/J_i$ is naturally isomorphic to $R[X_1,\ldots,X_{i-1},X_{i+1}\ldots,X_{r+1}]\cong R^{[r]}$. As $\bar\varphi(X_j)$ is a monomial for all $j\in \Lambda_\varphi\setminus\{i\}$, $|\Lambda_{\bar\varphi}|\ge d-1$.  Therefore, we can apply the induction hypothesis to conclude that there exist non-negative integers $m_0,n_0$, with $m_0\le r-(d-1)=r+1-d$ and $n_0\le d-1$, such that $\bar A$, as an $R$-algebra, is isomorphic to a $B_0$-retract of $B_0^{[m_0]}$, where $B_0\cong R^{[n_0]}$.\\
    Next, we assume that $\varphi(X_i)$ is a non-constant monomial for all $i\in \Lambda_\varphi$. If $i_0\in \Lambda_\varphi$, then we can write $\varphi(X_{i_0})$ as $\lambda M$, for some nonzero element $\lambda\in R$ and a monic monomial $M\in B$. If $X_i$ appears in $M$ for some $i$, then $\varphi(X_i)$ divides $\varphi(\lambda M)=\lambda M$. Since $R$ is an integral domain, every factor of a monomial is a monomial. Therefore, $\varphi(X_i)$ is a non-constant monomial for each $X_i$ that appears in $M$. If the degree of $\varphi(X_i)$ is more than one for some $X_i$ that appears in $M$, then the degree of $\varphi(\lambda M)$ would be greater than the degree of $\lambda M$, leading to a contradiction. Therefore, if $X_i$ appears in $M$, then $\varphi(X_i)$ is a linear monomial. Let $X_s$ be a variable that appears in $M$, so that $\varphi(X_s)=\mu X_t$ for some $t$ and a nonzero element $\mu\in R$. Since $\mu X_t=\varphi(\mu X_t)=\mu\varphi(X_t)$, it follows that $\varphi(X_t)=X_t$. Let $R_t:=R[X_t]\subseteq A$. Then $\varphi$ can be realized as an $R_t$-retraction of $R_t[X_1,\ldots,X_{t-1},X_{t+1},\ldots,X_{r+1}]\cong R_t^{[r]}$, say $\varphi_t$. Note that $\Lambda_\varphi\setminus\{t\}\subseteq\Lambda_{\varphi_t}$ (But equality may not hold, since we are now treating $X_t$ as a constant!). Therefore, $|\Lambda_{\varphi_t}|\ge d-1$, and we can apply the induction hypothesis to conclude that there exist nonnegative integers $m_0,n_0$, with $m_0\le r-(d-1)=r+1-d$ and $n_0\le d-1$, such that $A$, as an $R_t$-algebra (and hence, in particular, as an $R$-algebra), is isomorphic to a $B_0$-retract of $B_0^{[m_0]}$, where $B_0=R_t^{[n_0]}\cong R^{[n_0+1]}$.    
    
\end{proof}

Let $R$ be an integral domain and $n$ a non-negative integer. Let $B:=R[X_1,\ldots,X_n,Y_1,Y_2,Y_3]$ be the polynomial ring in $n+3$ variables over $R$, and $\varphi:B\longrightarrow B$ be an $R$-retraction with image $A$. If $\varphi(X_i)$ is zero or a monomial for all $i=1,\ldots,n$, then by  {\it proposition \ref{genprop2}}, there exist non-negative integers $m_0,n_0$, with $m_0\le 3$ and $n_0\le n$, such that $A$, as an $R$-algebra, is isomorphic to a $B_0$-retract of $B_0^{[m_0]}$, say $A_0$, where $B_0\cong R^{[n_0]}$. Clearly, $0\le \text{tr.deg}_{B_0}A_0\le m_0$. Note that $B_0$ is algebraically closed in $B_0^{[m_0]}$, and by \cite[Lemma 1.3]{costa1977retracts}, $A_0$ is algebraically closed in $B_0^{[m_0]}$. Therefore, if $\text{tr.deg}_{B_0}A_0=0$, then $A_0=B_0\cong R^{[n_0]}$; and if $\text{tr.deg}_{B_0}A_0=m_0$, then $A_0=B_0^{[m_0]}\cong R^{[m_0+n_0]}$. In either case, $A$ is isomorphic to a polynomial ring over $R$. So, we are left to investigate the cases where $1\le\text{tr.deg}_{B_0}A_0<m_0\le 3$, i.e., when $m_0=2$ and $\text{tr.deg}_{B_0}A_0=1$, and when $m_0=3$ and $\text{tr.deg}_{B_0}A_0=1$ or $2$. We consider these cases in the following corollary, where we tacitly use the fact that $A$ is isomorphic to $A_0$ as $R$-algebras.

\begin{corollary}
\label{gencor1}
With the same set-up and notation as in the above discussion

\begin{enumerate}
    \item if $R$ is a Noetherian domain and $\text{tr.deg}_{B_0}A_0=1$, then $A$ is an $\mathbb A^1$-fibration over $B_0$.
    \item if $R$ is a Noetherian semi-normal domain and $\text{tr.deg}_{B_0}A_0=1$, then $A\cong \text{Sym}_{B_0}(J)$ for some invertible ideal $J$ of $B_0$.
    \item if $R$ is a UFD and $\text{tr.deg}_{B_0}A_0=1$, then $A\cong B_0^{[1]}\cong R^{[n_0+1]}$.
    \item if $R$ is a Noetherian domain containing $\mathbb Q$, $m_0=3$ and $\text{tr.deg}_{B_0}A_0=2$, then $A$ is an $\mathbb A^2$-fibration over $B_0$.
\end{enumerate}
\end{corollary}

\begin{proof}
    If $R$ is a Noetherian (semi-normal) domain or a UFD, then so is $B_0$. In view of this, the assertions follow immediately from {\it theorems \ref{co1} and \ref{5.4}}.
\end{proof}
    
%\end{corollary}

\begin{remark}
\label{genrem1}
We continue with the same set-up and notation as in {\it corollary \ref{gencor1}} and the preceding discussion. We cannot improve (i), since there exist $R$-retractions of $R^{[2]}$, where $R$ is a Dedekind domain, with image $A$ such that $\text{tr.deg}_RA=1$, but $A\ncong R^{[1]}$ (see, for example, \cite[Example 4.7]{chakraborty2021some}). This also shows that in (iii), we cannot replace a UFD by even a Dedekind domain. \\ 
An $\mathbb{A}^1$-fibration over a Noetherian semi-normal domain $R$ is a symmetric algebra of an invertible ideal of $R$. We cannot drop the condition of semi-normality in (ii), as \cite[Example 4.1]{costa1977retracts} gives an example of a one-dimensional complete Noetherian local domain $R$, and an $R$-retraction of $R^{[2]}$ with an image $A$ such that $\text{tr.deg}_RA=1$, but $A$ is not isomorphic to the symmetric algebra of an invertible ideal of $R$.\\
In general, $\mathbb A^2$-fibrations over Noetherian UFDs are not trivial. For example, it is well-known that {\it the tangent bundle of the real sphere $S^2$ is not trivial, though the co-ordinate ring of $S^2$ is a UFD}. In the language of algebra, it is equivalent to the fact that {\it $R=\mathbb R[x,y,z]:=\mathbb R[X,Y,Z]/(X^2+Y^2+Z^2-1)$, where $x,y,z$ are the residue classes of the variables $X,Y,Z$, respectively, is a UFD, but the module of K{\"a}hler differentials of $R$, defined as $\Omega_{R/\mathbb R}:=\frac{R\oplus R\oplus R}{<(x,y,z)>}$, is a projective module of rank two over $R$ that is not free} (for a proof, see Propositions 8.3 and 8.4 of \cite{samuel_UFD_TIFR}). As $R^3\cong \Omega_{R/\mathbb R}\oplus R$ as $R$-modules, $A:=\text{Sym}_R(\Omega_{R/\mathbb R})$ is an $R$-retract of $R^{[3]}$. Since $\text{tr.deg}_RA=\text{rk}_R(\Omega_{R/\mathbb R})=2$, it follows from \cite[Theorem 5.9]{chakraborty2021some} that $A$ is an $\mathbb A^2$-fibration over $R$. However, $A\ncong R^{[2]}$, as $\Omega_{R/\mathbb R}$ is not a free $R$-module. Thus, (iv) appears to be the best possible result in general.\\
We do not attempt to generalize our results to the case where four or more variables map to arbitrary polynomials under $\varphi$. If $k$ is a field of positive characteristic, Gupta proved in \cite{gupta2014cancellation}, \cite{gupta2014zariski} that for all $n\ge 4$, there exist $k$-retracts of $k^{[n]}$ that are not polynomial rings over $k$; and if $k$ is a field of characteristic zero, then the question whether $k$-retracts of $k^{[n]}$ are polynomial rings over $k$ is still an open problem for all $n\ge 4$.
    
\end{remark}

\section{Binomial Retracts of Polynomial Rings}
\label{sec:binom}

In this section, we describe the structure of binomial retracts of polynomial rings. 

\begin{definition}
\label{binomialretract}
Let $R$ be a ring and $B:=R[X_1,\ldots,X_n]$ the polynomial ring in $n$ variables over $R$. An $R$-retraction $\varphi:B\longrightarrow B$ is called a {\it binomial $R$-retraction} if $\varphi(X_i)$ contains at most two terms for all $i=1,\ldots,n$, i.e., if $\varphi(X_i)$ is a constant, or a non-constant monomial or a binomial for all $i$. The image of a binomial $R$-retraction is called a {\it binomial $R$-retract}. \\ 
Note that monomial $R$-retractions are special cases of binomial $R$-retractions.\\
A binomial $R$-retraction $\varphi:B\longrightarrow B$ is called a {\it proper binomial $R$-retraction} if $\varphi(X_i)$ is a binomial for all $i$, i.e., if each $\varphi(X_i)$ is a sum of exactly two monomials.
    
\end{definition}

It is difficult to describe the structure of binomial retracts of polynomial rings over an arbitrary ring $R$. So, throughout this section, all rings that we consider are going to be integral domains.\\ First, we observe that the study of binomial retracts essentially boils down to the study of binomial retracts that arise from proper binomial retractions.

\begin{lemma}
    \label{binomlem0}
    Let $R$ be an integral domain and $B:=R[X_1,\ldots,X_n]$ the polynomial ring in $n$ variables over $R$. Let $\varphi:B\longrightarrow B$ be an $R$-retraction, not necessarily a binomial $R$-retraction. If $\varphi(X_i)$ is zero or a monomial for some $i$, then either $\varphi$ fixes $X_i$ or there exists a positive integer $j\in \{1,\ldots,n\}$ such that $X_j$ does not appear in $\varphi(X_j)$.
\end{lemma}

\begin{proof}
If $X_i$ does not appear in $\varphi(X_i)$, then we are done. Otherwise, we can write $\varphi(X_i)$ as $\varphi(X_i)=X_iM$ for some nonzero monomial $M\in B$. Then $X_iM=\varphi(X_iM)=X_iM\varphi(M)$, implying that $\varphi(M)=1$. If $M\in R$, then $M=\varphi(M)=1$, which implies that $\varphi$ fixes $X_i$. Otherwise, $M$ is divisible by a variable, say $X_j$. Then $\varphi(X_j)$ divides $\varphi(M)=1$, implying that $\varphi(X_j)$ is a unit of $R$. In particular, $X_j$ does not appear in $\varphi(X_j)$.    
\end{proof}

\begin{lemma}
    \label{binomlem1}
    Let $R$ be an integral domain, and $B:=R[X_1,\ldots,X_n]$ the polynomial ring in $n$ variables over $R$. Let $\varphi:B\longrightarrow B$ be a binomial $R$-retraction with image $A$. Then there exist non-negative integers $m_0,n_0$, with $m_0+n_0\le n$, and a proper binomial $B_0$-retraction $\tilde \varphi:B_0^{[m_0]}\longrightarrow B_0^{[m_0]}$, where $B_0\cong R^{[n_0]}$, such that $A$, as an $R$-algebra, is isomorphic to the image of $\tilde \varphi$. 
\end{lemma}

\begin{proof}
    We use induction on $n$, the number of variables. For $n=1$, since $R$ is an integral domain, it follows from {\it theorem \ref{co2}} that $\varphi(X)=X$ or $\varphi(X)\in R$. In the first case, we can take $M_0=0$ and $n_0=1$, and in the second case, we can take $m_0=n_0=0$. So, we assume that the assertion is true for some positive integer $r\ge 1$, and let $n=r+1$. 
    If $\varphi(X_i)$ is a binomial for all $i=1,\ldots,n$, then there is nothing to prove. Otherwise, there exists a positive integer $i$ such that $\varphi(X_i)$ is zero or a monomial. Then, by {\it lemma \ref{binomlem0}}, $\varphi$ fixes $X_i$ or there exists a positive integer $j$ such that $X_j$ does not appear in $\varphi(X_j)$.\\
    If $\varphi$ fixes $X_i$, let $R_i:=R[X_i]\subseteq A$, so that $B:=R_i[X_1,\ldots,X_{i-1},X_{i+1},\ldots,X_{r+1}]=R_i^{[r]}$. Clearly, we can realize $\varphi$ as a binomial $R_i$-retraction of $B$. Therefore, by the induction hypothesis, there exist non-negative integers $m_0,n_0$, with $m_0+n_0\le r$, and a proper binomial $B_0$-retraction $\tilde \varphi:B_0^{[m_0]}\longrightarrow B_0^{[m_0]}$, where $B_0\cong R_i^{[n_0]}=R^{[n_0+1]}$, such that $A$, as an $R_i$-algebra (and hence, as an $R$-algebra), is isomorphic to the image of $\tilde \varphi$. \\
    In the second case, let $X_j$ be a variable that does not appear in $\varphi(X_j)$. Let $J\subseteq B$ be the ideal generated by $X_j-\varphi(X_j)$. As $J\subseteq\text{ker }\varphi$, by {\it theorem \ref{5.4}\textcolor{blue}{(v)}}, $\varphi$ induces an $R$-retraction $\varphi':B/J\longrightarrow B/J$, such that the natural projection $B\longrightarrow B/J$ induces an $R$-algebra isomorphism between $A$ and $A':=\varphi'(B/J)$. Note that $B/J$ is naturally isomorphic to $B':=R[X_1,\ldots,X_{j-1},X_{j+1},\ldots,X_{r+1}]=R^{[r]}$, and $\varphi':B'\longrightarrow B'$ is a binomial $R$-retraction. Therefore, by the induction hypothesis, there exist non-negative integers $m_0,n_0$, with $m_0+n_0\le r$, and a proper binomial $B_0$-retraction $\tilde \varphi:B_0^{[m_0]}\longrightarrow B_0^{[m_0]}$, where $B_0\cong R^{[n_0]}$, such that $A'$, as an $R$-algebra, is isomorphic to the image of $\tilde \varphi$. Since $A\cong A'$ as $R$-algebras, the assertion follows.    
    %If $X_i$ does not appear in $\varphi(X_i)$, then we can consider the ideal $J\subseteq B$ generated by $X_i-\varphi(X_i)$. 
\end{proof}

The next couple of results help us to understand the structure of proper binomial $R$-retractions of a polynomial ring over $R$.

\begin{lemma}
    \label{binomlem2}
    Let $R$ be an integral domain, and $B:=R[X_1,\ldots,X_n]$ the polynomial ring in $n$ variables over $R$. Let $\varphi:B\longrightarrow B$ be an $R$-retraction with an image $A$. If $\text{ker }\varphi\cup A$ contains a non-constant monomial, then there exists a variable $X_j$ such that $\varphi(X_j)\in R$ or $\varphi(X_j)=X_j$. In particular, if $\varphi$ is a proper binomial $R$-retraction, then the following assertions hold.
    \begin{enumerate}
        \item The union $A\cup \text{ker }\varphi$ does not contain any non-constant monomial.
        \item The order of $\varphi(X_i)$ is at least one for all $i=1,\ldots,n$. In other words, $\varphi(B_+)\subseteq B_+$.\\
        In particular, if $M$ is a non-constant monomial of $B$, then $\text{ord }M\le \text{ord }\varphi(M)$, where the equality holds iff $\varphi(X_i)$ has order one for every variable $X_i$ that appears in $M$.
    \end{enumerate}
\end{lemma}

\begin{proof}    
    Let $M\in B$ be a non-constant monomial. Since $\text{ker }\varphi$ is a prime ideal of $B$, if $M\in\text{ ker } \varphi$ then $\text{ker }\varphi$ must contain at least one variable that appears in $M$, say $X_j$; and then $\varphi(X_j)=0$. Next, suppose that $M\in A$. Since $R$ is an integral domain, every factor of a monomial is also a monomial. Therefore, if $X_i$ is a variable that appears in $M$, then $\varphi(X_i)$ is a monomial. If none of the variables appearing in $M$ is mapped to a constant, then comparing the degrees of $M$ and $\varphi(M)$, we see that the image of every variable appearing in $M$ is a nonzero constant multiple of a variable. So, if $X_i$ appears in $M$ for some $i$, then $\varphi(X_i)=aX_j$ for some $j$ and a nonzero constant $a\in R$. Then $aX_j=\varphi(aX_j)=a\varphi(X_j)$, which implies that $\varphi(X_j)=X_j$.    
    \vspace{.4mm}\\
    (i) If $\varphi$ is a proper binomial retraction, then $\varphi(X_i)$ is a binomial for all $i=1,\ldots,n$. Therefore, the assertion immediately follows from the first part.
    \vspace{.4mm}\\
    (ii) If possible, let $\varphi(X_i)=c+\lambda M$ for some $i$, where $c,\lambda\in R$ are nonzero elements and $M$ is a non-constant monic monomial. Then $c+\lambda M=\varphi(X_i)=\varphi^2(X_i)=c+\lambda\varphi(M)$, implying that $M\in A$, which contradicts (i). \\
    The last assertion is obvious.
\end{proof}

\begin{lemma}
    \label{binomlem3}
    Let $R$ be an integral domain, and $B:=R[X_1,\ldots,X_n]$ the polynomial ring in $n$ variables over $R$. Let $\varphi:B\longrightarrow B$ be a proper binomial $R$-retraction. If $f:=\varphi(X_i)$ for some $i$, then the following assertions hold.
    \begin{enumerate}
        \item If $\text{ord }f=1$, then $f$ is a linear homogeneous polynomial.
        \item If $\text {ord }f\ge 2$, then $X_i$ does not appear in $\varphi(X_j)$ for all $j=1,\ldots,n$.\\
        In this case, $B_i:=R[X_1,\ldots,X_{i-1},X_{i+1},\ldots,X_n]$ is invariant under $\varphi$, and the induced proper binomial $R$-retraction $\varphi|_{B_i}:B_i\longrightarrow B_i$ has the same image as $\varphi$.
    \end{enumerate}
\end{lemma}

\begin{proof}
(i) If possible, suppose that $\text{ord }f=1$, but $f$ is not a linear homogeneous polynomial. Then $f$ can be written as $f=aX_k+bM$ for some $k\in \{1,\ldots,n\}$ and nonzero elements $a,b\in R$ such that $M$ is a monic monomial of degree $\ge 2$. Then $aX_k+bM=\varphi(X_i)=\varphi^2(X_i)=a\varphi(X_k)+b\varphi(M)$. By {\it lemma \ref{binomlem2}}\textcolor{blue}{(ii)}, the order of $\varphi(M)$ is at least two. So, comparing the linear terms of the two sides, we see that $\varphi(X_k)=X_k+b_kM_k$ for some nonzero element $b_k\in R$ and a non-constant monic monomial $M_k$. But then $X_k+b_kM_k=\varphi(X_k+b_kM_k)=X_k+b_kM_k+b_k\varphi(M_k)$, implying that $\varphi(M_k)=0$, which contradicts {\it lemma \ref{binomlem2}\textcolor{blue}{(i)}}. Therefore, if $\text{ord }f=1$, then $f$ is a linear homogeneous polynomial.\\
(ii) We assume that $\text{ord }f\ge 2$. If possible, suppose that $X_i$ appears in $\varphi(X_j)$ for some $j\in \{1,\ldots,n\}$. Let $\varphi(X_j)=M_1+M_2$, where $M_1,M_2$ are distinct non-constant monomials, not necessarily monic. Without loss of generality, we assume that $X_i$ divides $M_2$. If $d:=\text{ord }\varphi(X_j)$, then either $M_1$ or $M_2$ has order $d$. If $d=\text{ord }M_2<\text{ord }M_1$, then, by {\it lemma \ref{binomlem2}}\textcolor{blue}{(ii)}, both $\varphi(M_1)$ and $\varphi(M_2)$ have orders at least $d+1$. This is a contradiction, since $\varphi(X_j)=M_1+M_2=\varphi(M_1+M_2)=\varphi(M_1)+\varphi(M_2)$. So, we assume that $d=\text{ord }M_1\le \text{ord }M_2$. By (i) and {\it lemma \ref{binomlem2}}\textcolor{blue}{(ii)}, $\text{ord }M_i\le \varphi(M_i)$ for $i=1,2$, where equality holds iff every variable that appears in $M_i$ is mapped to a linear homogeneous binomial under $\varphi$. Since $X_i$ appears in $M_2$, $\text{ord }\varphi(M_2)>d$. Therefore, the order of $\varphi(M_1)$ is $d$, and it is a homogeneous polynomial with at least two terms. Since $\text{ord }\varphi(M_2)> \text{ord }\varphi(M_1)=\text{deg }\varphi(M_1)=d$, unless $\varphi(M_2)=0$, $\varphi(X_j)=\varphi(M_1)+\varphi(M_2)$ contains at least three terms, which is a contradiction. On the other hand, if $\varphi(M_2)=0$, then {\it lemma \ref{binomlem2}}\textcolor{blue}{(i)} leads to a contradiction.\\
The last assertion is obvious.
\end{proof}

We introduce a few definitions that will help  us in studying the structure of proper binomial $R$-retracts of a polynomial ring over $R$.

\begin{definition}
    Let $R$ be an integral domain and $B:=R[X_1,\ldots,X_n]$ the polynomial ring in $n$ variables over $R$. Let $\varphi:B\longrightarrow B$ be a proper binomial $R$-retraction. By {\it lemma \ref{binomlem3}\textcolor{blue}{(ii)}}, if the order of $\varphi(X_i)$ is at least $2$ for some $i$, then $X_i$ does not appear in the retract, and we can view $\varphi(B)$ as the image of the induced proper binomial $R$-retraction of $R[X_1,\ldots,X_{i-1},X_{i+1},\ldots,X_n]$. It allows us to ignore every variable that does not appear in the image of $\varphi$. So, using {\it lemma \ref{binomlem3}}, we implicitly assume here, as well as in the next {\it lemma}, that every variable appears in the retract, which implies that $\varphi(X_i)$ is a linear binomial for all $i=1,\ldots,n$. Then $\varphi$ induces an idempotent $R$-linear operator of the free $R$-module of rank $n$ generated by the set of variables $\mathfrak B:=\{X_1,\ldots,X_n\}$.\\
    
    We call a set $S\subseteq\mathfrak B$ {\it $\varphi$-invariant} if the $R$-subalgebra $R[X_i]_{i\in S}\subseteq B$ is invariant under $\varphi$, or equivalently, if, for all $X_i\in S$, $\varphi(X_i)$ is an $R$-linear combination of (two) elements of $S$. A $\varphi$-invariant set $S\subseteq \mathfrak B$ is called {\it minimal} if $S$ is minimal among all non-empty proper $\varphi$-invariant subsets of $\mathfrak B$ with respect to set inclusion.\\    
    If $\varphi(X_i)=aX_j+bX_k$ for some $i$, where $a,b\in R$ are nonzero elements, then we define $\mathcal V_i$ as the set of all variables that appear in $\varphi(X_j)$ or $\varphi(X_k)$. Since $\varphi$ is a binomial retraction and $aX_j+bX_k=a\varphi(X_j)+b\varphi(X_k)$, it follows that $|\mathcal V_i|\le 4$ and $\{j,k\}\subseteq \mathcal V_i$.
    
\end{definition}

\begin{lemma}
    \label{binomlem4}
    Let $R$ be an integral domain and $B:=R[X_1,\ldots,X_n]$ the polynomial ring in $n$ variables over $R$. Let $\mathfrak B:=\{X_1,\ldots,X_n\}$. If $\varphi:B\longrightarrow B$ is a proper binomial $R$-retraction, then the following assertions hold.
    
    \begin{enumerate} 
        \item Every non-empty $\varphi$-invariant subset of $\mathfrak B$ contains at least two elements; and a variable $X_i$ is contained in a $\varphi$-invariant set consisting of two elements iff $X_i$ appears in $\varphi(X_i)$.
        \item For each $i$, $\mathcal V_i$ is a $\varphi$-invariant set consisting of two or three elements, where the second possibility cannot occur, unless $\text{ch }R=2$.\\
        In particular, if $S\subseteq \mathfrak B$ is a minimal $\varphi$-invariant set, then $S=\mathcal V_i$ for all $X_i\in S$, and $S$ contains two or three elements, where the second possibility can occur only if $\text{ch }R=2$.
        \item For each $i=1,\ldots,n$, $\varphi(X_i)$ is contained in the $R$-linear span of a minimal $\varphi$-invariant subset of $\mathfrak B$, namely $\mathcal V_i$.
        \item If $S,T$ are minimal $\varphi$-invariant subsets of $\mathfrak B$, then $S=T$ or $S\cap T=\emptyset$.
        \item If $\text{ch }R=2$ and $S:=\{X_j,X_k,X_l\}$ is a $\varphi$-invariant subset of $\mathfrak B$ containing three elements, then the corresponding $R$-algebra $R':=R[X_j,X_k,X_l]$ is $\varphi$-invariant, and the induced proper binomial $R$-retraction $\varphi':=\varphi|_{R'}:R'\longrightarrow R'$ has an image that is isomorphic to $R^{[2]}$.
        
    \end{enumerate}
\end{lemma}

\begin{proof}
    (i) Since $\varphi(X_i)$ is an $R$-linear combination of two elements of $\mathfrak B$ for all $i=1,\ldots,n$, every non-empty $\varphi$-invariant subset of $\mathfrak B$ contains at least two elements. If $X_i$ is contained in a $\varphi$-invariant set consisting of two elements, say $\{X_i,X_j\}$, then it is clear from the definition of $\varphi$-invariant sets that both $\varphi(X_i)$ and $\varphi(X_j)$ are $R$-linear combinations of $X_i$ and $X_j$ with nonzero coefficients. In particular, $X_i$ appears in $\varphi(X_i)$, and $\mathcal V_i=\mathcal V_j=\{X_i,X_j\}$.\\ 
    Conversely, let $\varphi(X_i)=aX_i+bX_j$ for some $i$, where $a,b$ are nonzero elements of $R$. Then $aX_i+bX_j=a\varphi(X_i)+b\varphi(X_j)$, which implies that $b\varphi(X_j)=(1-a)(aX_i+bX_j)$. Since $\varphi(X_j)$ is an $R$-linear combination of two elements of $\mathfrak B$ and $\mathfrak B$ is linearly independent over $R$, it follows that $\varphi(X_j)$ is an $R$-linear combination of $X_i$ and $X_j$; and consequently, $\mathcal V_i=\mathcal V_j=\{X_i,X_j\}$ is $\varphi$-invariant.
    \vspace{.4mm}\\
    (ii) Let $\varphi(X_i)=aX_j+bX_k$ for some nonzero elements $a,b\in R$. If $X_i$ appears in $\varphi(X_i)$, i.e., if $j=i$ or $k=i$, then it follows from (i) that  $\mathcal V_i=\{j,k\}$ is a $\varphi$-invariant set containing two elements, and we are done.\\ 
    So, we assume $j,k\in \{1,\ldots,n\}\setminus\{i\}$. We know that $\{j,k\}\subseteq \mathcal V_i$ and $|\mathcal V_i|\le4$. If $|\mathcal V_i|=2$, then $\mathcal V_{i}=\{j,k\}$ and both $\varphi(X_j)$ and $\varphi (X_k)$ are nonzero $R$-linear combinations of $X_j$ and $X_k$. Therefore, $\mathcal V_i=\{j,k\}$ is a $\varphi$-invariant set consisting of two elements. On the other hand, if $|\mathcal V_i|=4$, then $\varphi(X_j)$ and $\varphi(X_k)$ do not share any common variable. But then $\varphi(X_i)=aX_j+bX_k=a\varphi(X_j)+b\varphi(X_k)$ contains four terms, which is a contradiction. Therefore, if $|\mathcal V_i|>2$, then $\mathcal V_i$ must contain exactly three elements.\\ 
    Let $\mathcal V_i=\{X_j,X_k,X_l\}$, where $l\in \{1,\ldots,n\}\setminus\{j,k\}$. If $X_j$ appears in $\varphi(X_j)$ and $X_k$ appears in $\varphi(X_k)$, then $\mathcal V_j=\{X_j,X_l\}$ and $\mathcal V_k=\{X_k,X_l\}$; for otherwise, if $X_l$ does not appear in both $\varphi(X_j)$ or $\varphi(X_k)$, but only in one of them, then $X_l$ would appear in $\varphi(X_i)$, leading to a contradiction. But then, by (i), $\mathcal V_l=\{j,l\}=\{k,l\}$, which is a contradiction. Therefore, without loss of generality, we are left with two cases, where $X_j$ appears in $\varphi(X_j)$ but $X_k$ does not appear in $\varphi(X_k)$, and where $X_j$ does not appear in $\varphi(X_j)$ and $X_k$ does not appear in $\varphi(X_k)$. In the first case, there exist nonzero elements $a_1,a_1,a_2,b_2\in R$ such that $\varphi(X_j)=a_1X_j+b_1X_k$ and $\varphi(X_k)=a_2X_j+b_2X_l$. But then \[a_1X_j+b_1X_k=\varphi^2(X_j)=a_1\varphi(X_j)+b_1\varphi(X_k)=a_1(a_1X_j+b_1X_k)+b_1(a_2X_j+b_2X_l),\]
    which is a contradiction, since $X_l$ does not appear in $\varphi(X_j)$.\\
    In the second case, there exist nonzero elements $c_1,d_1,c_2,d_2\in R$ such that $\varphi(X_j)=c_1X_k+d_1X_l$ and $\varphi(X_k)=c_2X_j+d_2X_l$. Then \[c_1X_k+d_1X_l=\varphi^2(X_j)=c_1\varphi(X_k)+d_1\varphi(X_l)=c_1(c_2X_j+d_2X_l)+d_1\varphi(X_l),\]
    which implies that $d_1\varphi(X_l)=c_1X_k-c_1c_2X_j+(d_1-c_1d_2)X_l$. Since $\varphi(X_l)$ is a binomial, it follows that $d_1=c_1d_2$ and $d_1\varphi(X_l)=c_1(X_k-c_2X_j)$. Therefore, $c_1$ divides $d_1$ and vice versa. Consequently, there exists a unit $\lambda\in R$ such that $c_1=\lambda d_1$, $d_2=\lambda^{-1}$ and $\varphi(X_l)=\lambda(X_k-c_2X_j)$. Therefore, $\mathcal V_i=\{X_j,X_k,X_l\}$ is a $\varphi$-invariant set containing three elements. Note that the binomial $R$-retraction $\varphi$ induces an idempotent $R$-linear operator of the free $R$-module $RX_j\oplus RX_k\oplus RX_l$ of rank $3$, say $\pi_\varphi$. Let $M_{\pi_\varphi}$ be the matrix representation of $\pi_\varphi$ with respect to the ordered basis $X_j,X_k,X_l$, Then $M_{\pi_\varphi}$, regarded as a $3\times 3$ matrix over $K:=Q(R)$, the field of fractions of $R$, is an idempotent matrix. Therefore, the rank of $M_{\pi_\varphi}$ is the same as its trace. But it is clear from the construction that $M_{\pi_\varphi}$ is a nonzero idempotent matrix of rank two whose trace is zero, which is possible only if $\text{ch }R=2$.\\
    Since $\mathcal V_i$ is a $\varphi$-invariant set for all $i=1,\ldots,n$, the second assertion follows immediately from the first one.
    \vspace{.4mm}\\
    (iii) This is obvious since $\mathcal V_i$ contains the two variables that appear in $\varphi(X_i)$. 
    \vspace{.4mm}\\
    (iv) Let $S,T$ be minimal $\varphi$-invariant subsets of $\mathfrak B$. If $X_i\in S\cap T$ for some $i$, then it follows from (ii) that $S=T=\mathcal V_i$.
    \vspace{.4mm}\\
    (v) Clearly, $R':=R[X_j,X_k,X_l]$ is invariant under $\varphi$. Since $S$ is a minimal $\varphi$-invariant set containing three elements, it follows from (i) that none of the variables in $S$ appears in its image under $\varphi$. Therefore, there exist nonzero elements $a_1,a_2,a_3,b_1,b_2,b_3$ of $R$ such that $\varphi(X_j)=a_1X_k+b_1X_l$, $\varphi(X_k)=a_2X_j+b_2X_l$ and $\varphi(X_l)=a_3X_j+b_3X_k$. A computation similar to the one in the last part of (ii) reveals that there exists a unit $\lambda\in R$ such that 
    \[\varphi(X_j)=b_1(\lambda X_k+X_l),\ \varphi(X_k)=a_2X_j+\lambda^{-1}X_l\ \text{ and }\ \varphi(X_l)=\lambda(X_k+a_2X_j).\]
    Since $\varphi=\varphi^2$, $a_2b_1=\lambda^{-1}$, which implies that both $a_2$ and $b_1$ are units of $R$. Clearly, the image of $\varphi':=\varphi|_{R'}$ is the symmetric $R$-algebra over the projective $R$-module $P:=R\varphi(X_j)+R\varphi(X_k)+R\varphi(X_l)$, which has rank two. As $\varphi(X_j)$ and $\varphi(X_k)$ are linearly independent over $R$ and $\varphi(X_l)=\varphi(X_j)/b_1+\lambda\varphi(X_k)$, it follows that $P=R\varphi(X_j)\oplus R\varphi(X_k)\cong R^{2}$; and consequently, the image of $\varphi'$ is isomorphic to $\text{Sym}_R(R^2)\cong R^{[2]}$.     
      
\end{proof}

\begin{remark}
    Let $R$ be an integral domain of characteristic $2$, and let $\varphi:R[X_1,X_2,X_3]\longrightarrow R[X_1,X_2,X_3]$ be a proper binomial $R$-retraction, defined as $\varphi(X_1):=X_2+X_3$, $\varphi(X_2):=X_1+X_3$ and $\varphi(X_3):=X_1+X_2$. Then $\mathfrak B:=\{X_1,X_2,X_3\}$ is already a minimal $\varphi$-invariant set. So, the case where $\text{ch }R=2$ indeed requires a separate treatment, as is done above.
\end{remark}

Now we describe the structure of proper binomial $R$-retracts of a polynomial ring over $R$.

\begin{theorem}
    \label{binomthm}
    Let $R$ be an integral domain, and $B:=R[X_1,\ldots,X_n]$ the polynomial ring in $n$ variables over $R$. Let $\varphi:B\longrightarrow B$ be a proper binomial $R$-retraction with an image $A$ such that $\text{tr.deg}_RA=r$. Then there exist a non-negative integer $s$ and projective $R$-modules $P_1,\ldots,P_s$ of rank one, each generated by two elements, such that 
    \begin{enumerate}
        \item if $\text{ch }R\neq 2$, then $A\cong \text{Sym}_R(P_1\oplus\cdots\oplus P_s)$ as $R$-algebras, where $s=r$.
        \item if $\text{ch }R=2$, then there exists a non-negative integer $t\ge s$ satisfying $2t=r+s$ such that $A$ is isomorphic to $\text{Sym}_R(P_1\oplus\cdots\oplus P_s)^{[2(t-s)]}$ as $R$-algebras.
    \end{enumerate}
    In particular, if $R$ is a UFD then irrespective of the characteristic of $R$, $A\cong R^{[r]}$ as $R$-algebras.    
    
\end{theorem}

\begin{proof}
    In view of {\it lemma \ref{binomlem3}}, we ignore all variables that do not appear in the image of $\varphi$, and assume that $\varphi(X_i)$ is a linear binomial for all $i$.
    Let $\mathfrak B:=\{X_1,\ldots,X_n\}$ and let $\Gamma\subseteq \mathfrak B$ be the union of minimal $\varphi$-invariant subsets of $\mathfrak B$. By {\it lemma \ref{binomlem4}}, $B_\Gamma:=R[X_i]_{i\in \Gamma}$ is $\varphi$-invariant and $\varphi(B)\subseteq B_\Gamma$. Since $\varphi^2=\varphi$, it follows that $A=\varphi(B_\Gamma)$. Therefore, replacing $B$ by $B_\Gamma$ and $\varphi$ by $\varphi|_{B_\Gamma}$, we assume that every variable of $B$ is contained in a minimal $\varphi$-invariant subset of $\mathfrak B$, the set of variables. Using {\it lemma \ref{binomlem4}}, we get a partition of $\mathfrak B$ into minimal $\varphi$-invariant sets $\Gamma_{1},\ldots,\Gamma_{s},\Gamma_{s+1},\ldots,\Gamma_{t}$ such that $|\Gamma_{i}|=2$ for all $1\le i\le s$ and $|\Gamma_{i}|=3$ for all $s+1\le i\le t$, where $s\le t$ are non-negative integers satisfying $2s+3(t-s)=n$, the number of variables. Note that if $\text{ch }R\neq 2$, then $s=t$. Let $B_i:=R[X_l]_{l\in \Gamma_i}$ for all $1\le i\le t$. Then each $B_i$ is $\varphi$-invariant and $B=B_1\otimes_R\cdots\otimes_RB_t$. For each $i\in \{1,\ldots,t\}$, let $\varphi_i:=\varphi|_{B_i}$ and $A_i:=\varphi_i(B_i)$. Then each $\varphi_i$ is a proper binomial $R$-retraction, and it is clear that $A\cong A_1\otimes_R\cdots\otimes_RA_t$. By {\it lemma \ref{binomlem4}\textcolor{blue}{(v)}}, $A_i\cong R^{[2]}$ for all $i=s+1,\ldots,t$. For a positive integer $i\le s$, let $\Gamma_i:=\{i_1,i_2\}$. Since $\varphi_i$ is a proper binomial $R$-retraction of $B_i=R[X_{i_1},X_{i_2}]$, both $\varphi_i(X_{i_1})$ and $\varphi_i(X_{i_2})$ are $R$-linear combinations of $X_{i_1}$ and $X_{i_2}$. Since $A_i\neq R$ and $A_i\neq B_i$, $\text{tr.deg}_RA_i=1$ for all $i=1,\ldots,s$. Note that $\varphi_i$ induces an $R$-linear endomorphism of the free $R$-module $RX_{i_1}\oplus RX_{i_2}\cong R^2$, say $\phi_i$, whose image is $P_i:=R\varphi(X_{i_1})+R\varphi(X_{i_2})$. It is clear that $A_i=\text{Sym}_R(P_i)$. Since $\phi_i$ is a non-trivial $R$-linear idempotent endomorphism, $R^2=\text{ker }\phi_i\oplus P_i$, implying that $P_i$ is a projective $R$-module of rank one that is generated by two elements. Therefore, by {\it proposition \ref{baki}}, 
    \[A\cong \text{Sym}_R(P_1)\otimes_R\cdots\otimes_R\text{Sym}_R(P_s)\otimes_RR^{[2(t-s)]}\cong \text{Sym}_R(P_1\oplus\cdots\oplus P_s)\otimes_R R^{[2(t-s)]},\]
    which is isomorphic to $\text{Sym}_R(P_1\oplus\cdots\oplus P_s)^{[2(t-s)]}$. Comparing the transcendence degrees of $A$ and $\text{Sym}_R(P_1\oplus\cdots\oplus P_s)^{[2(t-s)]}$ over $R$, we see that $r=2t-s$, or equivalently, $2t=r+s$. \\
    Finally, if $R$ is a UFD, then, by {\it theorem \ref{co1}}\textcolor{blue}{(ii)}, $A_i\cong R^{[1]}$ for all $i=1,\ldots,s$, which implies that $A\cong R^{[2t-s]}=R^{[r]}$ as $R$-algebras. Alternatively, since every finitely generated projective $R$-module of rank one over a UFD is free, we can conclude that in this case $P_i\cong R$ for all $i=1,\ldots,s$; and consequently, $A\cong \text{Sym}_R(R^s)\otimes_RR^{[2(t-s)]}\cong \text{Sym}_R(R^{2t-s})\cong R^{[r]}$.
    \end{proof}

\begin{remark}
   \label{binomrem}
    We continue with the same set-up and notation as in {\it theorem \ref{binomthm}}.
    \begin{enumerate}
        \item If $\text{ch }R\neq 2$, then $A\cong \text{Sym}_R(P_1\oplus\cdots\oplus P_r)$, where each $P_i$ is a projective $R$-module of rank one, generated by two elements; and if $\text{ch }R=2$, then $A\cong\text{Sym}_R(P_1\oplus\cdots\oplus P_s)^{[2(t-s)]}$, where $r=2t-s$. By {\it proposition \ref{baki}},
        \[\text{Sym}_R(P_1\oplus\cdots\oplus P_s)^{[2(t-s)]}\cong \text{Sym}_R(P_1\oplus\cdots\oplus P_s)\otimes_R\text{Sym}_R\big(R^{2(t-s)}\big),\]
        which is isomorphic to the symmetric algebra $\text{Sym}_R\big(P_1\oplus\cdots\oplus P_s\oplus R^{2(t-s)}\big)$. Therefore, even when $\text{ch }R=2$, the retract $A$ is isomorphic to the symmetric algebra of a projective $R$-module, which is a direct sum of $r$ projective modules of rank one, each generated by at most two elements.
        \item Since $\varphi(B_+)\subseteq B_+$, and the image of $\varphi$ is a graded $R$-subalgebra of $B$, we could have applied {\it theorem \ref{5.4}\textcolor{blue}{(iv)}} to conclude that the retract $A$ is isomorphic to $\text{Sym}_R(M) $ for some finitely generated projective $ R $-submodule $ M$ of the free $R$-module $\bigoplus_{i=1}^nRX_i$. However, proper binomial $R$-retractions are special where we could, in fact, show that $M$ is a direct sum of projective $R$-modules of rank one. Also, if $R$ is a UFD, then a projective $R$-module $M$, in general, may not be free; but we directly showed that, in this case, the retract $A$ is isomorphic to a polynomial ring over $R$. 
    \end{enumerate} 
\end{remark}

\begin{corollary}
    \label{binomcor}
    \label{binomcor}
    Let $R$ be an integral domain, and $B:=R[X_1,\ldots,X_n]$ the polynomial ring in $n$ variables over $R$. Let $\varphi:B\longrightarrow B$ be a binomial $R$-retraction with image $A$, and $m:=\text{tr.deg}_RA$. Then there exist non-negative integers $m_0,n_0$, with $m_0+n_0=m$, such that $A$, as an $R$-algebra, is isomorphic to the symmetric algebra $\text{Sym}_{R^{[n_0]}}(P_1\oplus\cdots\oplus P_{m_0})$, where each $P_i$ is a projective $R^{[n_0]}$-module of rank one that is generated by two elements. \\
    If $R$ is, moreover, 
    \begin{enumerate}
    \item a semi-normal domain, then $A\cong \text{Sym}_R(R^{n_0}\oplus Q_1\oplus\cdots\oplus Q_{m_0})$, where each $Q_i$ is a finitely generated projective $R$-module of rank one such that $P_i\cong R^{[n_0]}\otimes_RQ_i$.
    \item a UFD, then $A\cong R^{[m]}$ as $R$-algebras.
    \end{enumerate}
\end{corollary}

\begin{proof}
By {\it lemma \ref{binomlem1}}, there exist non-negative integers $n_0,n_1$, with $n_0+n_1\le n$, such that $A$, as an $R$-algebra, is isomorphic to the image of a proper binomial $R$-retraction $\varphi':B_0^{[n_1]}\longrightarrow B_0^{[n_1]}$, where $B_0\cong R^{[n_0]}$. Note that if $R$ is a UFD, then so is $B_0$. Let $A_0$ be the image of $\varphi'$ and $m_0:=\text{tr.deg}_{B_0}A_0$. Then, by {\it theorem \ref{binomthm}} and {\it remark \ref{binomrem}\textcolor{blue}{(i)}}, $A_0$ is isomorphic to the symmetric $B_0$-algebra $\text{Sym}_{B_0}(P_1\oplus\cdots\oplus P_{m_0})$, where each $P_i$ is a projective $B_0$-module of rank one that is generated by two elements; and if $B_0$ is a UFD, then $A_0\cong B_0^{[m_0]}=R^{[m_0+n_0]}$. Since $A\cong A_0$ as $R$-algebras, it follows that $m_0+n_0=m$.\\
Now assume that $R$ is semi-normal. By {\it theorem \ref{swan}}, each $P_i$ is extended from $R$, i.e., for all $i=1,\ldots,m_0$, there exists a projective $R$-module, say $Q_i$, of rank one, such that $P_i\cong B_0\otimes Q_i$. Therefore, $P_1\oplus\cdots\oplus P_{m_0}\cong B_0\otimes_R(Q_1\oplus\cdots\oplus Q_{m_0})$ as $B_0$-modules, and by {\it proposition \ref{baki}},%\textcolor{blue}{(i)}},
\[A_0\cong B_0\otimes_R\text{Sym}_R(Q_1\oplus\cdots\oplus Q_{m_0})\cong R^{[n_0]}\otimes_R\text{Sym}_R(Q_1\oplus\cdots\oplus Q_{m_0})\cong \text{Sym}_R(R^{n_0}\oplus Q_1\oplus\cdots\oplus Q_{m_0}).\]    
\end{proof}

\begin{remark}
We keep the set-up and notation of {\it corollary \ref{binomcor}}.
\begin{enumerate}
    \item Let $R:=\mathbb C[[t^2,t^3]]$, and $B:=R[X,Y,Z]$. Since
    \[\mathcal A:=\begin{bmatrix}
        1-t^4X^4 & t^2X^2+t^3X^3\\
        (t^2X^2-t^3X^3)(1+t^2X^2) & t^4X^4
    \end{bmatrix}\]
    is a non-trivial $2\times 2$ idempotent matrix over $R[X]$, we can define a binomial $R$-retraction $\varphi:B\longrightarrow B$ that fixes $X$, and sends $Y,Z$ to $(1-t^4X^4)Y+(t^2X^2+t^3X^3)Z$ and $\{(t^2X^2-t^3X^3)(1+t^2X^2)\}Y+t^4X^4Z$, respectively. Then $\varphi$ can be realized as a proper binomial $R[X]$-retraction of $\big(R[X]\big)[Y,Z]$. Therefore, the retract $A:=\varphi(B)$ is isomorphic to $\text{Sym}_{R[X]}(P_\mathcal A)$, where $P_\mathcal A$ is the projective $R[X]$-submodule of rank one of the free $R[X]$-module of rank two, given by the image of $\mathcal A$. Since $R$ is local, every finitely generated projective $R$-module is free. Therefore, if $A$ is isomorphic to the symmetric algebra of a projective $R$-module over $R$, then $A\cong R^{[2]}$. However, $A/(X)\cong R[Y]$ as $R$-algebras. Therefore, if $A\cong R^{[2]}$, then, by \cite[Theorem 3.9]{bhatwadekar}, $A= R[X]^{[1]}$, which is a contradiction, since $P_\mathcal A$ is not a free $R[X]$-module. Therefore, in {\it corollary \ref{binomcor}}, the retract $A$ may not be isomorphic to a symmetric algebra of a projective module over $R$, unless $R$ is semi-normal.
    \item Let $R$ be a Dedekind domain that is not a PID (For example, we can take $R:=\frac{\mathbb R[X,Y]}{(X^2+Y^2-1)}$, the co-ordinate ring of the real circle of radius one.). Then $R$ contains a nonzero ideal $J$ that is not principal. Since every ideal of $R$ is an invertible ideal, $J$ is a projective $R$-module of rank one. With every ideal of $R$ being generated by two elements, there exists an idempotent $R$-linear endomorphism $\pi:R^2\longrightarrow R^2$ whose image is isomorphic to $J$. By the universal property of symmetric algebras, $\pi$ naturally extends to a proper binomial $R$-retraction of $\text{Sym}_R(R^2)\cong R^{[2]}$ whose image is isomorphic to the symmetric $R$-algebra $\text{Sym}_R(J)$. Since $J$ is not a free $R$-module, it follows that $\text{Sym}_R(J)\ncong R^{[1]}$. Therefore, the results proved in {\it theorem \ref{binomthm}} and {\it corollary \ref{binomcor}} are perhaps the best possible results in general.
\end{enumerate} 
    
\end{remark}

\section{Open Problems}

\label{sec:open}

In this section, we record a couple of open problems, both inspired by the following example, which is a special case of an example that is due to Hamann (\cite[Example on Page 14]{Hamann}).

Let $\mathbb F_p$ be the finite field containing $p$ elements, where $p>0$ is a prime, and let $R:=\mathbb F_p+t^p\mathbb F_p[[t]]\subseteq\mathbb F_p[[t]]$. Let $R[U,V]$ be the polynomial ring in two variables over $R$. Let $\varphi:R[U,V]\longrightarrow R[U,V]$ be the $R$-retraction, defined by \[\varphi(U):=U-t^{p+1}V^p +t^{p+1}(U+tV)^{p^2}\ \text{ and }\ \varphi(V):=U^p +t^p(U+tV)^{p^2}.\] 
Then it is shown that the image of $\varphi$ is not isomorphic to the symmetric algebra of a projective $R$-module (of rank one).

In \cite[Corollary 3.8]{costa1977retracts} it is proved that {\it if $R$ is an Artinian ring, then every $R$-retract of $R^{[2]}$ is isomorphic to the symmetric algebra of a projective $R$-module of rank $\le 2$}. The above example shows that, in general, an $R$-retract of $R^{[2]}$ may not be isomorphic to the symmetric algebra of a projective module over $R$. However, we proved that {\it if $R$ is an integral domain, then every binomial $R$-retract of $R^{[2]}$ is isomorphic to the symmetric algebra of a projective $R$-module of rank $\le 2$} ({\it corollary \ref{diffcor}} and {\it theorem \ref{binomthm}}). Therefore, we ask the following question.

{\bf Question 1.} {\it If $R$ is a ring, is every binomial $R$-retract of $R^{[2]}$ isomorphic to the symmetric algebra of a projective $R$-module of rank $\le 2$?}

In view of {\it corollary \ref{diffcor}}, it is enough to consider proper binomial $R$-retractions; and using {\it lemma \ref{genelem}}, we can further assume that $R$ does not contain any non-trivial idempotent elements.

Now, suppose $R$ be an integral domain and $\varphi : R[X_1, X_2] \longrightarrow R[X_1, X_2]$ an $R$-retraction with image $A$. Let $t_1,t_2$ denote the number of terms in $\varphi(X_1)$ and $\varphi(X_2)$, respectively. Let $t_*:=\text{min }\{t_1,t_2\}$, $t^*:=\text{max }\{t_1,t_2\}$ and $\tilde t:=t_1+t_2$. We have seen that if $t_*\le 1$ or $\tilde t\le 4$, then $A$ is isomorphic to the symmetric algebra of a projective $R$-module of rank $\le 2$ ({\it corollary \ref{diffcor}}, {\it theorem \ref{binomthm}} and {\it remark \ref{binomrem}\textcolor{blue}{(i)}}). In particular, if $t^*\le 2$, i.e., if $\varphi$ is a binomial $R$-retraction, then $A$ is isomorphic to the symmetric algebra of a projective $R$-module of rank $\le 2$. Therefore, we ask the following question.

{\bf Question 2.} Let $R$ be an integral domain and $\varphi : R[X_1, X_2] \longrightarrow R[X_1, X_2]$ an $R$-retraction with image $A$. Let $t_1,t_2$ denote the number of terms in $\varphi(X_1)$ and $\varphi(X_2)$, respectively. Let $t_*:=\text{min}_i\{t_1,t_2\}$, $t^*:=\text{max}_i\{t_1,t_2\}$ and $\tilde t:=t_1+t_2$. Then what are the largest positive integers $c_*,c^*$ and $\tilde c$ such that $t_*\le c_*$, or $t^*\le c^*$ or $\tilde t\le \tilde c$ would imply that $A$ is isomorphic to the symmetric algebra of a projective $R$-module of rank $\le 2$?\\
And what if we further assume that $\text{ch }R=0$?

Hamann's example shows that for every prime $p$, there exists an integral domain $R$ of characteristic $p$ and an $R$-retraction of $R[X_1,X_2]$ such that $t_*=3,t^*=4$ and $\tilde t=7$, but the retract $A$ is not isomorphic to the symmetric algebra of a projective $R$-module. Therefore, for an integral domain $R$ with finite characteristic, we know that $1\le c_*\le 2$, $2\le c^*\le 3$ and $4\le \tilde c\le 6$. So, we ask whether these bounds can be further refined and whether there is any such example in characteristic zero.

\begin{Acknowledgement}
	The first author acknowledges Science and Engineering Research Board (SERB) for their MATRICS grant (Grant Number: [$\text{MTR}/2021/000297$]). The second author's research is funded by the National Board for Higher Mathematics (NBHM), Department of Atomic Energy (DAE), Government of India, Ref No: $ 0203/13(36)/2021-R \& D-II/13163 $.
\end{Acknowledgement}

\end{document}